\documentclass[11pt]{article}%
\usepackage{amsmath}
\usepackage{amsfonts}
\usepackage{amssymb}
\usepackage{graphicx}

\providecommand{\U}[1]{\protect\rule{.1in}{.1in}}

\numberwithin{equation}{section}
\newtheorem{theorem}{Theorem}[section]

\newtheorem{lemma}[theorem]{Lemma}

\newtheorem{proposition}[theorem]{Proposition}
\newtheorem{remark}[theorem]{Remark}

\def\<{\langle}
\def\>{\rangle}

\def\P{{\mathbb P}}
\def\Q{{\mathbb Q}}
\def\E{{\mathbb E}}

\def\R{{\mathbb R}}
\def\D{{\mathbb D}}
\def\N{{\mathbb N}}

\begin{document}

\parindent 0pt

\title{Regularity of Wiener functionals under\\
a H\"{o}rmander type condition of order one}
\author{ \textsc{Vlad Bally}\thanks{%
Universit\'e Paris-Est, LAMA (UMR CNRS, UPEMLV, UPEC), INRIA, F-77454
Marne-la-Vall\'ee, France. Email: \texttt{bally@univ-mlv.fr}.}
\smallskip \\
\textsc{Lucia Caramellino}\thanks{%
Dipartimento di Matematica, Universit\`a di Roma - Tor Vergata, Via della
Ricerca Scientifica 1, I-00133 Roma, Italy. Email: \texttt{%
caramell@mat.uniroma2.it}}\smallskip\\
}
\date{}
\maketitle

\parindent 0pt

{\textbf{Abstract.}}
We study the local existence and regularity of the density of the law of a functional on the Wiener space which satisfies a criterion that generalizes the H\"ormander condition of order one (that is, involving the first order Lie brackets) for diffusion processes.

\medskip

{\textbf{Keywords}}: Malliavin calculus; local integration by parts formulas; total variation distance; variance of the Brownian path.

\medskip

{\textbf{2010 MSC}}: 60H07, 60H30.


\section{Introduction}

H\"{o}rmander's theorem gives sufficient non degeneracy
assumptions under which the law of a diffusion process is absolutely
continuous with respect to the Lebesgue measure and has a smooth density. This
condition involves the coefficients of the diffusion process as well as the
Lie brackets up to an arbitrary order. The aim of this paper is to give a
partial generalization of this result to general functionals on the Wiener
space. We give in this framework a condition corresponding to the first
order H\"{o}rmander condition - we mean the condition which says that the
coefficients and the first Lie brackets span the space. Roughly speaking our regularity
criterion is as follows. Let $F$ be a functional on the Wiener space
associated to a Brownian motion $W=(W^{1},...,W^{d}).$ We denote by $D^{i}$
the Malliavin derivative with respect to $W^{i}$ and, for some $T>0$, we define%
\begin{equation}\label{1-intro}
\lambda(T)=\inf_{\left\vert \xi\right\vert =1}\Big(\sum_{i=1}^{d}\<
D_{T}^{i}F,\xi\> ^{2}+\sum_{i,j=1}^{d}\< D_{T}^{i}%
D_{T}^{j}F-D_{T}^{j}D_{T}^{i}F,\xi\> ^{2}\Big)%
\end{equation}
We fix $x$ and we suppose that
there exist $r,\lambda>0$ such that
\begin{equation}\label{ND-intro}
1_{\{\left\vert F-x\right\vert \leq r\}}(\lambda(T)-\lambda)\geq0\qquad \mbox{a.s.}
\end{equation}

Notice that, since $s\mapsto D_sF$ is defined as an element of $L^2([0,T])$,
the quantity $D_TF$ in (\ref{1-intro}) makes no sense. So, we will replace it by $\frac 1\delta\int_{T-\delta}^T\E_{T,\delta}(D_s F)ds$ for small values of $\delta$, where $\E_{T,\delta}$ denotes a suitable conditional expectation (see (\ref{NR}) for details).
Then, we actually replace (\ref{ND-intro}) with an asymptotic variant (see next Remark \ref{rem-Regularity}).

So, we assume that $F$ is five times
differentiable in Malliavin sense (actually in a slightly stronger sense) and
that the above non degeneracy condition holds for some $T>0.$ Then we prove that
the restriction of the law of $F$ to $B_{r/2}(x)$ is absolutely continuous and
has a smooth density.

The analysis of the Malliavin covariance matrix under the non degeneracy hypothesis \eqref{ND-intro} is based on an estimate concerning the variance of the Brownian path. This is done by using its Laplace transform, which has been studied by Donati-Martin and Yor \cite{[DY]}.  We employ also another important argument, which is  the regularity criterion for the law of a random variable given in \cite{[BC2]}: it allows one to use integration by parts formulas in an ``asymptotic way''.

The main result is Theorem \ref{Regularity}, and Section \ref{sect-main} is devoted to its proof, for which we use results on the variance of the Brownian path which are postponed to Appendix \ref{app-1}. In Section \ref{sect-use} we illustrate the result with an example from diffusion processes with coefficients which may depend on the path of the process.

At our knowledge there are not many results concerning general vectors on the
Wiener space - except of course the celebrated criterion given by Malliavin
and the Bouleau Hirsh criterion for the absolute continuity. Another
criterion proved by Kusuoka in \cite{[K]} and further generalized by Nourdin and
Poly \cite{[NP]} and Nualart, Nourdin and Poly \cite{[NNP]} concerns vectors living
in a finite number of chaoses. All these criterions suppose that the
determinant of the Malliavin covariance matrix is non null in a more or less
strong sense - but give no hint about the possible analysis of this condition.
This remains to be checked using ad hoc methods in each particular example. So
the main progress in our paper is to give a rather general condition under
which the above mentioned determinant behaves well.

\bigskip

\textbf{Acknowledgments.} We are grateful to E. Pardoux who made a remark which allowed
us to improve a previous version of our result.

\section{Existence and smoothness of the local density}\label{sect-main}

Let us recall some notations from Malliavin calculus (we refer to Nualart \cite{[N]} or Ikeda and Watanabe
\cite{[IW]}). We work on a probability space $(\Omega,\mathcal{F},P)$ with a $d$
dimensional Brownian motion $W=(W^{1},...,W^{d})$ and we denote by
$\mathcal{F}_{t}$ the standard filtration associated to $W.$ We fix a time-horizon $T_0>0$ and we denote by
$\D^{k,p}$ the space of the functionals on the Wiener space which are $k$ times
differentiable in $L^{p}$ in Malliavin sense on the time interval $[0,T_0]$ and we put $\D^{k,\infty}%
=\cap_{p\geq1}\D^{k,p}.$ For a multi index $\alpha=(\alpha_{1},...,\alpha
_{k})\in\{1,...,d\}^{k}$ and a functional $F\in \D^{k,p}$ we denote $D^{\alpha
}F=(D_{s_{1},...,s_{k}}^{\alpha}F)_{s_{1},...,s_{k}\in [0,T_0]}$ with
$D_{s_{1},...,s_{k}}^{\alpha}F=D_{s_{k}}^{\alpha_{k}}...D_{s_{1}}^{\alpha_{1}%
}F$. Moreover, for $|\alpha|=k$ we define the norms%
\begin{align}
\left\vert D^{\alpha}F\right\vert _{L^{p}[0,T_0]^{k}}^{p}  &
:=\int_{[0,T_0]^{k}}\left\vert D_{s_{1},...,s_{k}}^{\alpha}F\right\vert
^{p}ds_{1},...,ds_{k}\quad \mbox{and}\quad\label{M1}\\
\left\Vert F\right\Vert _{k,p}  &  =\left\Vert F\right\Vert _{p}%
+\sum_{r=1}^k\sum_{\vert \alpha\vert = r}\E(\left\vert D^{\alpha
}F\right\vert _{L^{2}[0,T_0]^{r}}^{p})^{1/p}.\nonumber
\end{align}
If $F=(F^{1},...,F^{n})$, we set
$$
\left\vert D^{\alpha}F\right\vert
_{L^{p}[0,T_0]^{k}}^{p}=\sum_{i=1}^{n}\left\vert D^{\alpha}F^{i}\right\vert
_{L^{p}[0,T_0]^{k}}^{p}\quad\mbox{and}\quad\left\Vert F\right\Vert _{k,p}=\sum_{i=1}%
^{n}\left\Vert F^{i}\right\Vert _{k,p}.
$$
Moreover we will use the following
seminorms:%
\begin{align*}
\left\vert \left\vert \left\vert F\right\vert \right\vert \right\vert
_{k,p,q}  &  =\sum_{r=3}^{k}\sum_{\left\vert \alpha\right\vert =r}\E(\left\vert
D^{\alpha}F\right\vert _{L^{q}[0,T_0]^{r}}^{p})^{1/p}\\
&  =\sum_{r=3}^{k}\sum_{\left\vert \alpha\right\vert =r}\E\Big(\Big(\int
_{[0,T_0]^{r}}\left\vert D_{s_{1},...,s_{k}}^{\alpha}F\right\vert
^{q}ds_{1}...ds_{r}\Big)^{p/q}\Big)^{1/p}.
\end{align*}
Notice that $\||\cdot|\|_{k,p,q}$ does not take into account $\left\Vert F\right\Vert
_{p}$ and the norm of the first two derivatives. Moreover, for $q=2$ we find out the usual norms but if $q>2$ the control given by $\left\vert \left\vert
\left\vert F\right\vert \right\vert \right\vert _{k,q,p}$ (on the derivatives
of order larger or equal to three) is stronger than the one given by
$\left\vert \left\vert F\right\vert \right\vert _{k,p}.$ We define the spaces%
\[
\D^{k,p}=\{F:\left\Vert F\right\Vert _{k,p}<\infty\},\qquad \D^{k,p,q}%
=\D^{k,p}\cap\{F:\left\Vert |F|\right\Vert _{k,p,q}<\infty\}.
\]

Clearly $\D^{k,p,q}\subset \D^{k,p}$ for $q>2$\ and for $q=2$ we have equality.
We also denote
\begin{equation}
\D^{k,\infty}=\cap_{p\geq1}\D^{k,p},\qquad \D^{k,\infty,q}=\cap_{p\geq1}%
\D^{k,p,q}, \qquad \D^{k,\infty,\infty}=\cap_{p\geq1}%
\cap_{q\geq2}\D^{k,p,q} \label{M2}%
\end{equation}

For $s<t$ we denote
\[
\mathcal{F}_{s}^{t}=\mathcal{F}_{s}\vee\sigma(W_{u}-W_{t},u\geq t)=\sigma
(W_{v},v\leq s)\vee\sigma(W_{u}-W_{t},u\geq t).
\]
Now, for a fixed instant $T\in (0,T_0]$, we denote by $\E_{T,\delta}$ the conditional expectation with respect to
$\mathcal{F}_{T-\delta}^{T}$ that is%
\begin{equation}\label{NR}
\E_{T,\delta}(\Theta)=\E(\Theta\mid\mathcal{F}_{T-\delta}^{T}).
\end{equation}
We will use the following slight extension of the Clark-Ocone formula:
for $F\in \D^{1,2}$ and for $0\leq\delta<T$ one has%
\begin{equation}
F=\E_{T,\delta}(F)+\sum_{i=1}^{d}\int_{T-\delta}^{T}\E_{T,T-s}(D_{s}^{i}%
F)dW_{s}^{i}. \label{M2'}%
\end{equation}
\eqref{M2'} is immediate for simple functionals, and then can be straightforwardly generalized to functionals in $\D^{1,2}$.

For $\delta\in (0,T)$, we consider a family of random vectors
$$
a(T,\delta
)=(a_{i}(T,\delta),a_{k,j}(T,\delta))_{i,k,j=1,...,d}
$$
and we assume that $a(T,\delta)$ is $\mathcal{F}_{T-\delta}^{T}$ measurable.
We denote%
\begin{equation}\label{R1'}
\begin{split}
&[a]_{i,j}(T,\delta)    =a_{i,j}(T,\delta)-a_{j,i}(T,\delta
),\\
&\overline{a}(T,\delta)  =\Big(\sum_{i=1}^{d}\left\vert a_{i}%
(T,\delta)\right\vert ^{2}+\sum_{i,j=1}^{d}\left\vert a_{i,j}(T,\delta
)\right\vert ^{2}\Big)^{1/2},\\
&\lambda(T,\delta)    =\inf_{\left\vert \xi\right\vert =1}\Big(\sum_{i=1}%
^{d}\left\langle a_{i}(T,\delta),\xi\right\rangle ^{2}+\sum_{i,j=1}%
^{d}\left\langle [a]_{i,j}(T,\delta),\xi\right\rangle ^{2}\Big).
\end{split}
\end{equation}
For $p\geq1,\alpha>0,0<\delta<T$ we define%
\begin{equation}\label{R1}
\begin{split}\varepsilon_{\alpha,p,T,\delta}(a,F) :=&\sum_{i=1}^{d}\Big(\E\Big(\Big(\frac{1}{\delta
}\int_{T-\delta}^{T}\left\vert \frac{\E_{T,\delta}(D_{s}^{i}F)-a_{i}(T,\delta
)}{\delta^{\frac{1}{2}+\alpha}}\right\vert ^{2}ds\Big)^{p}\Big)\Big)^{1/2p}+\\
&  +\sum_{i,j=1}^{d}\Big(\frac{1}{\delta^{2}}\int_{T-\delta}^{T}\int_{T-\delta
}^{s_1}\E\Big(\left\vert \frac{\E_{T,\delta}(D_{s_{2}}^{j}D_{s_{1}}^{i}F)-a_{i,j}%
(T,\delta)}{\delta^{\alpha/2}}\right\vert ^{2p}\Big)ds_{1}ds_{2}\Big)^{1/2p}.
\end{split}
\end{equation}

Our main result is the following.

\begin{theorem}\label{Regularity}
Let $F=(F^{1},..,F^{n})$ be $\mathcal{F}_{T_0}$-measurable with $F^{i}\in \D^{2,\infty},i=1,...,n.$ We fix $y\in \R^{n}$ and $r>0$ and we suppose that there exists
$\alpha,\lambda_{\ast}>0,\gamma\in\lbrack0,\frac{1}{2}),$ a time $T\in(0,T_0]$ and a
family $a(T,\delta)=(a_{i}(T,\delta),a_{i,j}(T,\delta))_{i,j=1,...,d},$ of
$\mathcal{F}_{T-\delta}^{T}$ measurable vectors\ such that for every $p\geq1$
\begin{equation}\label{R3}
\begin{split}
i) &\qquad\limsup_{\delta\rightarrow0}\varepsilon_{\alpha,p,T,\delta
}(a,F)  <\infty,\\
ii) & \qquad\limsup_{\delta\rightarrow0}\delta^{p\gamma}\E(\overline
{a}^{p}(T,\delta)) <\infty,\\
iii) & \qquad\limsup_{\delta\rightarrow0}\delta^{-p}\P(\{\left\vert
F-y\right\vert  \leq r\}\cap\{\lambda(T,\delta)<\lambda_{\ast}%
\})<\infty.
\end{split}
\end{equation}
Then the following statements hold.

\medskip

\textbf{A}. Suppose that $F^{i}\in \cup_{p>6}\D^{5,\infty,p},i=1,...,n.$ Then the law
of $F$ on $B_{r/2}(y):=\{x:\left\vert x-y\right\vert <r/2\}$ is absolutely
continuous with respect to the Lebesgue measure. We denote by $p_{F}$ the
density of the law.

\smallskip

\textbf{B}. Suppose that for some $k\geq 5$ one has $F^{i}\in \D^{k,\infty,\infty},i=1,...,n.$ Then
$$
p_{F}\in\cap_{p\geq1}W^{k-5,p}(B_{r/2}(y)).
$$
\end{theorem}

\begin{remark}\label{rem-Regularity}
Morally, $D^i_TF\sim \frac 1\delta\int_{T-\delta}^T\E_{T,\delta}(D^i_s F)ds$. Then, condition $i)$ in (\ref{R3}) says that we may replace $D^i_TF$ by $a_i(T,\delta)$, and we have a precise control of the error. The same for $D^i_T D^j_t F$, which is replaced by $a_{i,j}(T,\delta)$. Then, $iii)$ in (\ref{R3}) gives the asymptotic non-degeneracy condition in terms of $\lambda(T,\delta)$, which is associated to $a(T,\delta)$.
\end{remark}

\begin{remark}\label{rem-Regularity2}
Notice that we may ask the non-degeneracy condition $iii)$ in (\ref{R3}) to hold in any intermediary time $T\in(0,T_0]$ and not only for $T=T_0$ (we thank to E. Pardoux for a remark in this sense).

\end{remark}

The proof is postponed to Section \ref{proof}. We first need to state some preliminary results.

\subsection{A short discussion on the proof of Theorem \ref{Regularity}}

Let us give the main ideas and the strategy we are going to use to prove Theorem \ref{Regularity}.

We will look to the law of $F$ under $\P_{U}$ where $U$ is a localization random variable for the set $\{\left\vert F-y\right\vert \leq
r\}.$ We want to prove that this law is absolutely continuous with respect to
the Lebesgue measure - this implies that the law of $F$ restricted to
$\{\left\vert F-y\right\vert \leq r\}$ is absolutely continuous (and this is
our aim). In order to do it we proceed as follows: for each $\delta>0$ we
construct some localization random variables $U_{\delta}$ in such a way that
on the set $\{U_{\delta}\neq0\}$ the random variable $F$ has nice properties -
this means that we may control the Malliavin derivatives and the Malliavin
covariance matrix of $F$ on the set $\{U_{\delta}\neq0\}$. This allows us to build integration by parts formulas
for $F$ under $\P_{U_{\delta}}.$ The $L^{p}$ norms of the weights which appear
in these integration by parts formulas blow up as $\delta\rightarrow0$ but we
have a sufficiently precise control of the rate of the blow up. On the other
hand we will estimate the total variation distance between the law of $F$
under $\P_{U}$ and under $\P_{U_{\delta}}.$ We prove that this distance goes to
zero as $\delta\rightarrow0$ and we obtain sufficiently precise estimates of
the rate of convergence. Then we use Theorem 2.13 in \cite{[BC2]}, that we recall here in next Theorem \ref{regularity-new}, which
guarantees that if one may achieve a good equilibrium between the rate of the
blow up and the rate of convergence to zero, then one obtains a density for the
limit law.

It worth to stress that the strategy employed here is slightly different from
the usual one. In fact, in next (\ref{R4}) we decompose $F$ as $F=\E_{T,\delta}(F)+Z_\delta(a)+R_\delta$ and   one would expect that we approximate $F$ by $\E_{T,\delta}(F)+Z_\delta(a)$. But we do not proceed in this way. We keep all the time the
same random variable $F$ (which includes $R_\delta$) but we change the probability measure under which we
work in order to have a good localization: we replace $\P_{U}$ by
$\P_{U_{\delta}}.$ The decomposition $F=\E_{T,\delta}(F)+Z_\delta(a)+R_{\delta}$ is not used in order to produce the approximation $\E_{T,\delta}(F)+Z_\delta(a)$ but just to analyze the properties for $F$ itself under different localizations given in $\P_{U_{\delta}}.$ As we will see soon, such a  decomposition appears as a Taylor expansion of order one in which $Z_\delta(a)$ represents the principal term and $R_\delta$ is a reminder in the sense that it is small on the set $\{U_\delta\neq 0\}$.

\subsection{Preliminary results}
Let $F\in\D^{4,2}$. Using twice Clark Ocone formula (\ref{M2'}), we obtain
\begin{equation}
F-\E_{T,\delta}(F)=Z_{\delta}(a)+R_{\delta}(F) \label{R4}%
\end{equation}
with%
\begin{equation}
Z_{\delta}(a)=\sum_{i=1}^{d}a_{i}(T,\delta)(W_{T}^{i}-W_{T-\delta}^{i}%
)+\sum_{i,j=1}^{d}a_{i,j}(T,\delta)\int_{T-\delta}^{T}(W_{s}^{i}-W_{T-\delta
}^{i})dW_{s}^{j} \label{R6}%
\end{equation}
and $R_{\delta}(F)=R_{\delta}^{\prime}(F)+R_{\delta}^{\prime\prime}(F)$ with
\begin{equation}\label{R5}
\begin{split}
R_{\delta}^{\prime}(F) = &  \sum_{i=1}^{d}\int_{T-\delta}^{T}(\E_{T,\delta
}(D_{s}^{i}F)-a_{i}(T,\delta))dW_{s}^{i}\\
&  +\sum_{i,j=1}^{d}\int_{T-\delta}^{T}\int_{T-\delta}^{s_{1}}(\E_{T,\delta
}(D_{s_{2}}^{j}D_{s_{1}}^{i}F)-a_{i,j}(T,\delta))dW_{s_{2}}^{j}dW_{s_{1}}%
^{i}\\
R_{\delta}^{\prime\prime}(F) = &  \sum_{i,j,k=1}^{d}\int_{T-\delta}^{T}%
\int_{T-\delta}^{s_{1}}\int_{T-\delta}^{s_{2}}\E_{T,T-s_{3}}(D_{s_{3}}%
^{k}D_{s_{2}}^{j}D_{s_{1}}^{i}F)dW_{s_3}^{k}dW_{s_{2}}^{j}dW_{s_{1}}%
^{i}
\end{split}
\end{equation}

Since $T$ and $\delta$ are fixed we will use in the following shorter notation
\[
a_{i}=a_{i}(T,\delta),\qquad a_{i,j}=a_{i,j}(T,\delta),\qquad\overline
{a}=\overline{a}(T,\delta).
\]
We will use the Malliavin calculus restricted to $W_{s},$ $s\in\lbrack
T-\delta,T].$ Straightforward computations give
\begin{equation}\label{R5'}
D_{s}^{j}Z_{\delta}(a)    =a_{j}+\sum_{i\neq j}[a]_{i,j}(W_{s}^{i}
-W_{T-\delta}^{i})+r_{j}, \mbox{ with }
r_{j}    =\sum_{i=1}^{d}a_{i,j}(W_{T}^{i}-W_{T-\delta}^{i}).
\end{equation}

We denote
\begin{equation} \label{R5'''}
\begin{array}{c}
\displaystyle
q_{1}(W)=\left\vert W_{T}-W_{T-\delta}\right\vert ,\qquad q_{2}(W)=\frac
{1}{\delta}\int_{T-\delta}^{T}\left\vert W_{s}-W_{T-\delta}\right\vert
^{2}ds,\smallskip\\
\displaystyle
G_{\delta}=\int_{T-\delta}^{T}\left\vert D_{s}R_{\delta
}\right\vert ^{2}ds%
\end{array}
\end{equation}
and we define
\begin{equation}\label{Lambda}
\Lambda_{T,\delta}=\Big\{q_{1}(W)\leq\frac{1}{8\overline{a}}
\frac{\sqrt{\lambda_\ast}}d\Big\}\cap\Big\{q_{2}(W)\leq1\Big\}\cap\Big\{G_{\delta}\leq
\frac{\lambda_{\ast}}{34}\delta^{2}\Big\}\cap\Big\{\lambda(T,\delta)\geq\lambda_{\ast
}\Big\}
\end{equation}
We set $\sigma_{F,T,\delta}$ as the Malliavin covariance matrix of $F$
associated to the Malliavin derivatives restricted to $W_{s},$ $s\in\lbrack
T-\delta,T]$ that is%
\begin{equation}
\sigma_{F}^{i,j}=\sigma_{F,T,\delta}^{i,j}=\int_{T-\delta}^{T}\left\langle
D_{s}F^{i},D_{s}F^{j}\right\rangle ds, \quad i,j=1,\ldots,n. \label{R5''}%
\end{equation}
The main step of the proof is the following estimate. It is based on an analysis of the variance of the Brownian path, which is done in Appendix \ref{app-1}.

\begin{lemma}\label{Determinant}
Let $F=(F^{1},...,F^{n})$ with $F^{i}\in \D^{4,2}$. Let $0\leq\delta<T$ be fixed and $\E_{T,\delta}$ be defined in \eqref{NR}. Then for every $p\geq1$
\begin{equation}
\E_{T,\delta}(1_{\Lambda_{T,\delta}}(\det\sigma_{F,T,\delta})^{-p})\leq
\frac{C_{n,p}}{\lambda_{\ast}^{pn}\delta^{2pn}} \label{R8a}%
\end{equation}
with%
\[
C_{n,p}=2\Gamma(p)\int_{\R^n}\left\vert \xi\right\vert ^{n(2p-1)}e^{-\frac{1}{34}\left\vert
\xi\right\vert ^{2}}d\xi.
\]

\end{lemma}

\textbf{Proof.} By using Lemma 7-29, pg 92 in \cite{[BGJ]}, for every $n\times n$
dimensional and non negative defined matrix $\sigma$ one has%
\[
(\det\sigma)^{-p}\leq\Gamma(p)\int\left\vert \xi\right\vert ^{n(2p-1)}%
e^{-\left\langle \sigma\xi,\xi\right\rangle }d\xi,
\]
so that
\[
\E_{T,\delta}((\det\sigma_{F})^{-p}1_{\Lambda_{T,\delta}})\leq\Gamma(p)\int\left\vert
\xi\right\vert ^{n(2p-1)}\E_{T,\delta}(1_{\Lambda_{T,\delta}}e^{-\left\langle \sigma
_{F}\xi,\xi\right\rangle })d\xi.
\]
Since $\Lambda_{T,\delta}\subset\{G_{\delta}\leq\frac{\lambda_{\ast}}{34}\delta^{2}\}$
we have%
\begin{align*}
\left\langle \sigma_{F}\xi,\xi\right\rangle  &  \geq\frac{1}{2}\left\langle
\sigma_{Z_{\delta}(a)}\xi,\xi\right\rangle -\left\langle G_{\delta}\xi
,\xi\right\rangle \geq\frac{1}{2}\left\langle \sigma_{Z_{\delta}(a)}\xi
,\xi\right\rangle -G_{\delta}\left\vert \xi\right\vert ^{2}\\
&  \geq\frac{1}{2}\left\langle \sigma_{Z_{\delta}(a)}\xi,\xi\right\rangle
-\frac{\lambda_{\ast}}{34}\delta^{2}\left\vert \xi\right\vert ^{2}%
\end{align*}
so that%
\[
\E_{T,\delta}((\det\sigma_{F})^{-p}1_{\Lambda_{T,\delta}})\leq\Gamma(p)\int\left\vert
\xi\right\vert ^{n(2p-1)}e^{\frac{\lambda_{\ast}}{34}\delta^{2}\left\vert
\xi\right\vert ^{2}}\E_{T,\delta}(1_{\Lambda_{T,\delta}}e^{-\left\langle \sigma
_{Z_{\delta}(a)}\xi,\xi\right\rangle })d\xi.
\]
We fix $\xi\in \R^n$ and we choose $j=j(\xi)$ such that
\[
\left\langle a_{j},\xi\right\rangle ^{2}+\sum_{i\neq j}\left\langle
[a]_{i,j},\xi\right\rangle ^{2}\geq\frac{\lambda_{\ast}}{d}\left\vert
\xi\right\vert ^{2}.
\]
This is possible because we are on the set $\Lambda_{T,\delta}\subset\{\lambda(T,\delta)\geq\lambda_{\ast}\}.$ Then by (\ref{R5'})
\begin{align*}
\left\langle \sigma_{Z_{\delta}(a)}\xi,\xi\right\rangle
&=\int_{T-\delta
}^{T}\left\langle D_{s}^{j}Z_{\delta}(a),\xi\right\rangle ^{2}ds\\
&=\int
_{T-\delta}^{T}\Big(  \left\langle a_{j},\xi\right\rangle +\left\langle
r_{j},\xi\right\rangle +\sum_{i\neq j}\left\langle [a]_{i,j},\xi\right\rangle
(W_{s}^{i}-W_{T-\delta}^{i})\Big)  ^{2}ds.
\end{align*}
We define
$$
\beta_j^{2}(\xi)=\sum_{i\neq j}\left\langle [a]_{i,j},\xi\right\rangle ^{2}
$$
and for $\beta_j^{2}(\xi)>0$,
\[
b_s(j,\xi)=\frac{1}{\beta_j(\xi)}\sum_{i\neq j}\left\langle [a]_{i,j}%
,\xi\right\rangle (W_{T-\delta+s}^{i}-W_{T-\delta}^{i}).
\]
Notice that $b(j,\xi)$ is a Brownian motion under $\P_{T,\delta}.$ We also set
$b_s(j,\xi)=0$ in the case $\beta_j^{2}(\xi)=0$. Then the
previous equality reads%
\[
\left\langle \sigma_{Z_{\delta}(a)}\xi,\xi\right\rangle
=\int_{0}^{\delta}\left(  \left\langle a_{j},\xi\right\rangle +\left\langle r_{j}%
,\xi\right\rangle +\beta_j(\xi)b_{s}(j,\xi)\right)  ^{2}ds.
\]
We use now Lemma \ref{Varience} in Appendix \ref{app-1} with $\alpha=\left\langle a_{j}%
,\xi\right\rangle ,\beta=\beta_j(\xi),r=\left\langle r_{j},\xi\right\rangle $
and $b_{s}=b_{s}(j,\xi).$ We have to check that the assumptions there are
verified. Using Cauchy-Schwarz inequality we obtain%
\begin{align*}
\frac{1}{\delta}\Big|\int_{0}^{\delta} b_{s}(j,\xi)ds\Big|
&\leq
\Big(\frac{1}{\delta}\int_{0}^{\delta}\left\vert b_{s}(\xi)\right\vert
^{2}ds\Big)^{1/2}\\
&\leq\Big(
\frac{1}{\delta}\int_{0}^{\delta}\left\vert W_{T-\delta
+s}-W_{T-\delta}\right\vert ^{2}ds
\Big)^{1/2}=\sqrt{q_{2}(W)}\leq1.
\end{align*}
Moreover, since $\alpha^{2}+\beta^{2}\geq\frac{\lambda_{\ast}}{d}\left\vert
\xi\right\vert ^{2}$ we have
$$
\left\vert r\right\vert ^{2}\leq\left\vert r_{j}\right\vert ^{2}\left\vert
\xi\right\vert ^{2}\leq d\overline{a}^{2}q_{1}^{2}(W)\left\vert \xi\right\vert
^{2}\leq\frac{1}{64}\frac{\lambda_{\ast}}d\left\vert \xi\right\vert ^{2}\leq\frac
{1}{64}(\alpha^{2}+\beta^{2}).
$$
So the hypothesis are verified: by using (\ref{A2}) we obtain%
\[
\E_{T,\delta}(1_{\Lambda_{T,\delta}}e^{-\left\langle \sigma_{Z_{\delta}(a)}\xi
,\xi\right\rangle })\leq2e^{-\frac{\delta^{2}}{17}(\left\vert \alpha
\right\vert ^{2}+\left\vert \beta\right\vert ^{2})}\leq2e^{-\frac{\delta
^{2}\lambda_{\ast}}{17d}\left\vert \xi\right\vert ^{2}}.
\]
We come back and we obtain
\begin{align*}
\E_{T,\delta}((\det\sigma_{F})^{-p}1_{\Lambda_{T,\delta}})
&\leq2\Gamma(p)\int\left\vert
\xi\right\vert ^{n(2p-1)}e^{\frac{\lambda_{\ast}}{34d}\delta^{2}\left\vert
\xi\right\vert ^{2}}e^{-\frac{\delta^{2}\lambda_{\ast}}{17d}\left\vert
\xi\right\vert ^{2}}d\xi\\
&=2\Gamma(p)\int\left\vert \xi\right\vert
^{n(2p-1)}e^{-\frac{\delta^{2}\lambda_{\ast}}{34}\left\vert \xi\right\vert
^{2}}d\xi\\
&=\frac{C_{n,p}}{\lambda_{\ast}^{pn}\delta^{2pn}}
\end{align*}
where the last equality easily follows by a change of variable .
$\square$

\medskip

We also need the following estimate.

\begin{lemma}\label{Norms}
Suppose that (\ref{R3}) $i)$ holds and let $G_\delta$ be defined as in \eqref{R5'''}.

\medskip

\textbf{A.} If $F^{i}\in \cup_{p>6}  \D^{4,\infty,p},i=1,...,n$, there exists $\varepsilon>0$ such that
\begin{equation}
\limsup_{\delta\rightarrow0}\,\delta^{-\varepsilon}\,\P(G_{\delta}%
\geq\delta^{2})<\infty. \label{R10}%
\end{equation}

\medskip

\textbf{B.} If $F^{i}\in \D^{4,\infty,\infty},i=1,...,n$ then (\ref{R10}) holds
for every $\varepsilon>0.$
\end{lemma}

\textbf{Proof}. \textbf{A.} Let $F\in (\D^{4,\infty,p})^n$ for some $p>6$. We recall that $R'_\delta$ and $R''_\delta$ are defined in \eqref{R5}.
We write $R_{\delta}^{\prime}%
(F)=\sum_{i=1}^{d}r_{\delta}^{i}+\sum_{i,j=1}^{d}r_{\delta}^{i,j}$ and
$R''_{\delta}=\sum_{i,j,k=1}^d r^{i,j,k}_\delta$,
with%
\begin{align*}
r_{\delta}^{i}
&  =\int_{T-\delta}^{T}(\E_{T,\delta}(D_{s}^{i}F)-a_{i}%
(T,\delta))dW_{s}^{i},\\
r_{\delta}^{i,j}
&  =\int_{T-\delta}^{T}\int_{T-\delta}^{s_{1}}(\E_{T,\delta
}(D_{s_{2}}^{j}D_{s_{1}}^{i}F)-a_{i,j}(T,\delta))dW_{s_{2}}^{j}dW_{s_{1}}^{i},\\
r^{i,j,k}_\delta
&=\int_{T-\delta}^{T}%
\int_{T-\delta}^{s_{1}}\int_{T-\delta}^{s_{2}}\E_{T,T-s_{3}}(D_{s_{3}}%
^{k}D_{s_{2}}^{j}D_{s_{1}}^{i}F)dW_{s_3}^{k}dW_{s_{2}}^{j}dW_{s_{1}}%
^{i}.
\end{align*}

\textbf{Step 1}. We estimate $G_{\delta}^{i}=\int_{T-\delta}^{T}\left\vert
D_{s}^{i}r_{\delta}^{i}\right\vert ^{2}ds.$ For $s\in\lbrack T-\delta,T]$\ we
have $D_{s}^{i}r_{\delta}^{i}=\E_{T,\delta}(D_{s}^{i}F)-a_{i}(T,\delta)$ so%
\[
G_{\delta}^{i}=\int_{T-\delta}^{T}\left\vert \E_{T,\delta}(D_{s}^{i}%
F)-a_{i}(T,\delta)\right\vert ^{2}ds.
\]
It follows that
\begin{align*}
\frac{1}{\delta^{\varepsilon}}\P(G^i_{\delta}    \geq\delta^{2})
&\leq\frac
{1}{\delta^{\varepsilon}}\delta^{-2p}\left\Vert G_{\delta}^{i}\right\Vert
_{p}^{p}
=\frac{1}{\delta^{\varepsilon}}\E\Big(\Big|\delta^{-1}\int_{T-\delta}%
^{T}\left\vert \frac{\E_{T,\delta}(D_{s}^{i}F)-a_{i}(T,\delta)}{\delta^{1/2}%
}\right\vert ^{2}ds\Big|^{p}\Big)\\
&  \leq\frac{1}{\delta^{\varepsilon}}\times\delta^{2\alpha p}\varepsilon
_{\alpha,p,T,\delta}^{2p}(a,F)
\end{align*}
and consequently, by our hypothesis (\ref{R3}) $i),$ this term satisfies
(\ref{R10}) for every $\varepsilon>0$ (it suffices to take $p$ sufficiently large).

\textbf{Step 2}. We estimate $G_{\delta}^{i,j}=\sum_{\ell=1}^d\int_{T-\delta}^{T}\left\vert
D_{s}^{\ell}r_{\delta}^{i,j}\right\vert ^{2}ds$. We have
\begin{align*}
D_{s}^{\ell}r_{\delta}^{i,j}    =
&1_{i=\ell}\int_{T-\delta}^{s}(\E_{T,\delta
}(D_{s_{2}}^{j}D_{s}^{i}F)-a_{i,j}(T,\delta))dW_{s_{2}}^{j}+\\
&+1_{j=\ell}
\int_{s}^{T}(\E_{T,\delta}(D_{s}^{p}D_{s_{1}}^{i}F)-a_{i,p}(T,\delta
))dW_{s_{1}}^{i}=:1_{i=\ell}u_{s}^{i,j}+1_{j=\ell}v_{s}^{i,j}.
\end{align*}
We have%
\begin{align*}
\E\Big(\Big|\int_{T-\delta}^{T}\left\vert u_{s}^{i,j}\right\vert ^{2}ds\Big|^{p}\Big)  &
\leq\delta^{p-1}\int_{T-\delta}^{T}\E\big(\left\vert u_{s}^{i,j}\right\vert
^{2p}\big)ds\\
&  \leq C\delta^{p-1}\int_{T-\delta}^{T}\E\Big(\Big|\int_{T-\delta}^{s}(\E_{T,\delta
}(D_{s_{2}}^{j}D_{s}^{i}F)-a_{i,j}(T,\delta))^{2}ds_{2}\Big|^{p}\Big)ds\\
&  \leq C\delta^{2p-2}\int_{T-\delta}^{T}\int_{T-\delta}^{s}\E\big(\big|\E_{T,\delta
}(D_{s_{2}}^{j}D_{s}^{i}F)-a_{i,j}(T,\delta)\big|^{2p}\big)ds_{2}ds\\
&  =C\delta^{2p+\alpha p}\frac{1}{\delta^{2}}\int_{T-\delta}^{T}\int
_{T-\delta}^{T}\E\Big(\left\vert \frac{\E_{T,\delta}(D_{s_{2}}^{j}D_{s_{1}}%
^{i}F)-a_{i,j}(T,\delta)}{\delta^{\alpha/2}}\right\vert ^{2p}\Big)ds_{1}ds_{2}\\
&  \leq C\delta^{2p+\alpha p}\varepsilon_{\alpha,p,T,\delta}^{2p}(a,F).
\end{align*}
Using Chebyshev inequality we obtain
\[
\P\Big(\int_{T-\delta}^{T}\left\vert u_{s}^{i,j}\right\vert ^{2}ds\geq\delta
^{2}\Big)\leq C\delta^{-2p}\delta^{2p+\alpha p}\varepsilon_{\alpha,p,T,\delta
}^{2p}(a,F)=C\delta^{\alpha p}\varepsilon_{\alpha,p,T,\delta}^{2p}(a,F).
\]
which by (\ref{R3}) $i),$ satisfies (\ref{R10}) for every $\varepsilon>0$. For
$v_{s}^{i,j}$ the argument is the same.

\textbf{Step 3}. We estimate $G_{\delta}^{i,j,k}=\sum_{\ell=1}^d\int_{T-\delta}^{T}\vert
D_{s}^{\ell}r_{\delta}^{i,j,k}\vert ^{2}ds$.
We have
\begin{align*}
D^\ell_sr^{i,j,k}_\delta
=&1_{i=\ell}\int_{T-\delta}^s\int_{T-\delta}^{s_2}
\E_{T,T-s_{3}}(D_{s_{3}}%
^{k}D_{s_{2}}^{j}D_{s}^{i}F)dW_{s_3}^{k}dW_{s_{2}}^{j}+\\
&+1_{j=\ell}
\int_{T-\delta}^{T}%
\int_{T-\delta}^{s_{2}}1_{s<s_1}
\E_{T,T-s_{3}}(D_{s_{3}}^{k}D_{s}^{j}D_{s_{1}}^{i}F)dW_{s_3}^{k}dW_{s_{1}}^{i}+\\
&+1_{k=\ell}
\int_{T-\delta}^{T}%
\int_{T-\delta}^{s_{1}}1_{s<s_2}\E_{T,T-s}(D_{s}%
^{k}D_{s_{2}}^{j}D_{s_{1}}^{i}F)dW_{s_{2}}^{j}dW_{s_{1}}%
^{i}+\\
&+\int_{T-\delta}^{T}%
\int_{T-\delta}^{s_{1}}\int_{T-\delta}^{s_{2}}1_{s<s_3}D^{\ell}_s\E_{T,T-s_{3}}(D_{s_{3}}%
^{k}D_{s_{2}}^{j}D_{s_{1}}^{i}F)dW_{s_3}^{k}dW_{s_{2}}^{j}dW_{s_{1}}%
^{i}\\
=:&1_{i=\ell}u_s^{i,j,k}+1_{j=\ell}v_s^{i,j,k}+1_{k=\ell}w_s^{i,j,k}+z_s^{i,j,k,\ell}.
\end{align*}
By using H\"{o}lder and Burkholder inequality as in step 1, one obtains
\begin{align*}
\E\Big(\Big|\int_{T-\delta}^{T}\left\vert u_{s}^{i,j,k}\right\vert ^{2}ds\Big|^{p}\Big)  &\leq
\delta^{3p-3}\int_{T-\delta}^T\int_{T-\delta}^T\int_{T-\delta}^T
\E\big(|D_{s_{3}}%
^{k}D_{s_{2}}^{j}D_{s}^{i}F|^{2p}\big)ds_3ds_{2}ds\\
&\leq \delta^{3p-3}|||F|||_{3,2p,2p}^{2p}.
\end{align*}
An identical bound holds for $\E(|\int_{T-\delta}^{T}\vert v_{s}^{i,j,k}\vert ^{2}ds|^{p})$ and $\E(|\int_{T-\delta}^{T}\vert w_{s}^{i,j,k}\vert ^{2}ds|^{p})$. As for $z^{i,j,k,\ell}$, one more further integral appears, so we get $\E(|\int_{T-\delta}^{T}\vert z_{s}^{i,j,k,\ell}\vert ^{2}ds|^{p})
\leq\delta^{4p-4}|||F|||_{4,2p,2p}^{2p}$. By summarizing, we get
$$
\E\Big(\Big|\int_{T-\delta}^{T}\left\vert D_sR''_\delta\right\vert ^{2}ds\Big|^{p}\Big)\leq \delta^{3p-3}|||F|||_{4,2p,2p}^{2p}
$$
so that for every $p>1$%
\[
\P\Big(\int_{T-\delta}^{T}\left\vert D_{s}^{\ell}R_{\delta}^{\prime\prime}\right\vert
^{2}ds\geq\delta^{2}\Big)\leq C\delta^{-2p}\delta^{3p-3}|||
F||| _{4,2p,2p}^{2p}=C\delta^{p-3}|||F||| _{4,2p,2p}^{2p}.
\]
Suppose first that $F^{i}\in\cup_{p>6}\D^{4,\infty,p}.$ Then we may find $p>3$
such that $||| F|||_{4,2p,2p}<\infty$
and consequently the above quantity is upper bounded by $C\delta^{p-3}.$
This means that (\ref{R10}) holds for $\varepsilon<p-3.$ If $F^{i}\in
\D^{4,\infty,\infty}$ then we may take $p$ arbitrary large and so we obtain
(\ref{R10}) for every $\varepsilon>0.$

$\square$

\medskip

We will also need the following property for $G_\delta$.

\begin{lemma}\label{Norms1bis}
If $F\in\D^{k+1,2p}$ then
$$
\|G_\delta\|_{k,p}
\leq C\big(\|F\|^2_{k+1,2p}+\delta\|\overline{a}(T,\delta)\|_{4p}^2\big),
$$
where $C$ denotes a constant depending on $k,p,d$ only.

\end{lemma}

\textbf{Proof.}
For $G\in (\D^{k,p})^n$, we set $|D^{(k)}G|=\sum_{\ell=0}^k\sum_{|\gamma|=\ell}|D^{\gamma}G|^2$, where, for $|\gamma|=\ell$,
$$
|D^{\gamma}G|^2=\int_{[0,T]^\ell}
|D^\gamma_{s_1\ldots s_\ell}G|^2ds_1\cdots ds_\ell,
$$
that is $|D^{\gamma}G|$ is the one given in \eqref{M2} with $p=2$. Here, the case $|\gamma|=0$, that is $\gamma=\emptyset$, reduces to the original random variable: $D^\emptyset G=G$ and $|D^{(0)}G|=|G|$.

In the following, we let $C$ denote a positive constant, independent of $\delta$ and the random variables we are going to write. And we let $C$ vary from line to line.

We take $G_\delta=\int_{t-\delta}^T|D_sR_\delta|^2ds$ and we first prove the following (deterministic) estimate: there exists a constant $C$ depending on $k$ and $d$ such that
\begin{equation}\label{Gdelta-est}
|D^{(k)}G_\delta|\leq C|D^{(k+1)}R_\delta|^2.
\end{equation}
For $k=0$, this is trivial. Consider $k=1$. One has $
D^i_uG_\delta
=\sum_{\ell=1}^d\int_{t-\delta}^T2D^\ell_sR_\delta D^i_uD^\ell_sR_\delta ds$, so that, by using the Cauchy–Schwarz inequality, we get
\begin{align*}
|DG_\delta|^2
&\leq 4\sum_{i,\ell=1}^d
\int_{T-\delta}^T\big|\int_{T-\delta}^T2D^\ell_sR_\delta D^i_uD^\ell_sR_\delta ds\Big|^2du\\
&\leq 4\sum_{i,\ell=1}^d
\int_{T-\delta}^T du
\int_{T-\delta}^T2|D^\ell_sR_\delta|^2ds
\int_{T-\delta}^T2|D^i_uD^\ell_sR_\delta|^2 ds\\
&
\leq
C|D^{(1)}R_\delta|^2|D^{(2)}R_\delta|^2\leq
C|D^{(2)}R_\delta|^4
\end{align*}
and \eqref{Gdelta-est} holds for $k=1$. For $k\geq 2$, we use the following straightforward formula: if $\alpha$ denotes a multi-index of length $k$, then
$$
D^\alpha G_\delta
=\sum_{\ell=1}^d\int_{T-\delta}^T
\Big(
2D^{\ell}_sR_\delta D^{\alpha}R_\delta
+\sum_{\beta\in\mathcal{P}_\alpha}
D^\beta D^{\ell}_sR_\delta D^{\alpha\setminus\beta}D^\ell_s R_\delta\Big)ds,
$$
where $\mathcal{P}_\alpha$ is the set of the non empty multi indeces $\beta$ which are a subset of $\alpha$ and $\alpha\setminus\beta$ stands for the multi index of length $|\alpha|-|\beta|$ given by eliminating from $\alpha$ the entries of $\beta$. By using the above formula and the Cauchy–Schwarz inequality, one easily  gets
\begin{align*}
\int_{[T-\delta,T]^k}\!\!|D^{\alpha}_{s_1,\ldots,s_k}G_\delta|^2ds_1\cdots ds_k
&\leq
C\Big(|D^{(1)}R_\delta|^2|D^{(k)}R_\delta|^2
+\!\sum_{r=1}^k |D^{(r+1)}R_\delta|^2|D^{(k-r+1)}R_\delta|^2\Big)\\
&\leq C |D^{(k+1)}R_\delta|^4
\end{align*}
and \eqref{Gdelta-est} follows. Passing to expectation in \eqref{Gdelta-est}, it follows that
$$
\|G_\delta\|_{k,p}\leq C \|R_\delta\|_{k+1,2p}^2
$$
and by recalling that $R_\delta=F-\E_{T,\delta}(F)-Z_\delta(a)$, we obtain
$$
\|G_\delta\|_{k,p}
\leq C\big(\|F\|^2_{k+1,2p}+\|Z_\delta(a)\|_{k+1,2p}^2\big).
$$
From \eqref{R6}, by using H\"older's inequality we get
\begin{align*}
\|Z_{\delta}(a)\|_{k+1,2p}
\leq
&
\sum_{i=1}^{d}\|a_{i}(T,\delta)\|_{4p}\|W_{T}^{i}-W_{T-\delta}^{i}\|_{k+1,4p}\\
&+\sum_{i,j=1}^{d}\|a_{i,j}(T,\delta)\|_{4p}\Big\|\int_{T-\delta}^{T}(W_{s}^{i}-W_{T-\delta
}^{i})dW_{s}^{j}\Big\|_{k+1,4p}\\
&\leq C \|\overline{a}(T,\delta)\|_{4p}\delta^{1/2}
+C \|\overline{a}(T,\delta)\|_{4p}\delta\\
&\leq
C\delta^{1/2}\,\|\overline{a}(T,\delta)\|_{4p},
\end{align*}
and the statement follows.
$\square$

\begin{remark}\label{Norms1}
If hypothesis \eqref{R3} $ii)$ holds then $\limsup_{\delta\to 0}\delta\|\overline{a}(T,\delta)\|_{4p}^2=0$  because in this case one takes $\gamma<1/2$, so that for $F\in(\D^{k+1,2p})^n$ one has
$$
\sup_{\delta>0}\|G_\delta\|_{k,p}<\infty.
$$
\end{remark}

\subsection{Localization}

We will use a localization argument from \cite{[BC2]} that we recall here$.$ We
consider a random variable $U$ taking values in $[0,1]$ and we denote
\begin{equation}
d\P_{U}=Ud\P.\label{loc1}%
\end{equation}
This is a non negative measure (but generally not a probability measure - one
must divide with $\E(U)$ to get a probability measure). We denote%
\begin{align}
&\left\Vert F\right\Vert _{U,p} :=\E_{U}(\left\vert F\right\vert
^{p})^{1/p}=\E(\left\vert F\right\vert ^{p}U)^{1/p}\quad \mbox{and}\quad
\label{loc2}\\
&\left\Vert F\right\Vert _{U,k,p} :=\left\Vert F\right\Vert _{U,p}%
+\sum_{r=1}^{k}\sum_{\left\vert \alpha\right\vert =r}\E_{U}(\left\vert
D^{\alpha}F\right\vert _{L^{2}[0,T_0]^{r}}^{p})^{1/p}.\nonumber
\end{align}
Clearly $\left\Vert F\right\Vert _{U,k,p}\leq\left\Vert F\right\Vert _{k,p}.$
For a random variable $F\in(\D^{1,2})^{n}$ we denote
\begin{equation}
\sigma_{U,F}(p)=\E_{U}((\det\sigma_{F})^{-p})^{1/p}.\label{loc3}%
\end{equation}
We assume that $U\in \D^{1,\infty}$ and that for every $p\geq1$%
\begin{equation}
m_{p}(U):=\E_{U}(\left\vert D\ln U\right\vert ^{p})<\infty.\label{loc4}%
\end{equation}
In Lemma 2.1 in \cite{[BC1]} we have proved the following:

\begin{lemma}
\label{IP}Assume that (\ref{loc4}) holds. Let $F\in(\D^{2,\infty})^{n}$ be such
that $\det\sigma_{F}\neq0$ on the set $\{U\neq0\}$. We denote $\widehat{\sigma
}_{F}$ the inverse of $\sigma_{F}$ and we assume that $\sigma_{U,F}%
(p)<\infty$ for every $p\in \N$. Then for every $V\in \D^{1,\infty}$ and every
$f\in C_{b}^{\infty}(\R^{n})$ one has%
\begin{equation}
\E_{U}(\partial_{i}f(F)V)=\E_{U}(f(F)H_{i,U}(F,V)) \label{loc5}%
\end{equation}
with%
\begin{equation}
H_{i,U}(F,V)=\sum_{j=1}^{n}\big(V\widehat{\sigma}_{F}^{j,i}LF^{j}-\left\langle
D(V\widehat{\sigma}_{F}^{j,i}),DF^{j}\right\rangle -V\widehat{\sigma}%
_{F}^{j,i}\left\langle D(\ln U),DF^{j}\right\rangle \big). \label{loc6}%
\end{equation}
Suppose that $\ln U\in \D^{k,\infty}.$\ Iterating (\ref{loc5}) one obtains for
a multi index $\alpha=(\alpha_{1},...,\alpha_{k})\in\{1,...,n\}^{k}$%
\begin{equation}
\E_{U}(\partial_{\alpha}f(F)V)=\E_{U}(f(F)H_{\alpha,U}(F,V)),
\mbox{ with $H_{\alpha,U}(F,V)=H_{\alpha_{k},U}(F,H_{\underline{\alpha}
,U}(F,V))$
}, \label{loc7}%
\end{equation}
where $\underline{\alpha}=(\alpha_{1},...,\alpha_{k-1})$.

\end{lemma}

We will use this result with a localization random variable $U$ constructed in
the following way. For $a\in(0,1)$ we define $\psi_{a}:\R_{+}\rightarrow \R_{+}$
by
\begin{equation}
\psi_{a}(x)=1_{[0,a)}(x)+1_{[a,2a)}(x)\exp\Big(1-\frac{a^{2}}{a^{2}-(x-a)^{2}}\Big).
\label{loc8}%
\end{equation}
Then for every multi index $\alpha$ and every $p\in \N$ there exists a
universal constant $C_{\alpha,p}$ such that%
\begin{equation}
\sup_{x\in \R_{+}}\psi_{a}(x)\left\vert \partial_{\alpha}\ln\psi_{a}%
(x)\right\vert ^{p}\leq\frac{C_{\alpha,p}}{a^{p\left\vert \alpha\right\vert }%
}. \label{loc9}%
\end{equation}
Let $a_{i}>0$ and $Q_{i}\in \D^{1,p},i=1,...,l$ and $U=\prod_{i=1}^{l}%
\psi_{a_{i}}(Q_{i}).$ As an easy consequence of (\ref{loc9}) we obtain the
following estimates
\begin{equation}
m_{p}(U)\leq C\sum_{i=1}^{l}\frac{1}{a_{i}^{p}}\left\Vert Q_{i}\right\Vert
_{U,1,p}^{p}\leq C\sum_{i=1}^{l}\frac{1}{a_{i}^{p}}\left\Vert Q_{i}\right\Vert
_{1,p}^{p} \label{loc10}%
\end{equation}
where $C$ is a universal constant. And moreover, for every $k,p\in \N$ there
exists a universal constant $C$ such that%
\begin{equation}
\left\Vert \ln U\right\Vert _{U,k,p}\leq C\sum_{i=1}^{l}\frac{1}{a_{i}^{k}%
}\left\Vert Q_{i}\right\Vert _{k,p}. \label{loc11}%
\end{equation}

The function $\psi_{a}$ is suited for localization around zero. In order to
localize far from zero we have to use the following alternative version:
\begin{equation}
\phi_{a}(x)=1_{[a,\infty)}(x)+1_{[a/2,a)}(x)\exp\Big(1-\frac{a^{2}}{(2x-a)^{2}}\Big).
\label{loc8'}%
\end{equation}
The property (\ref{loc9}) holds for $\phi_{a}$ as well. And if one employs both
$\psi_{a_{i}}$ and $\phi_{a_{i}}$ in the construction of $U$, that is if one sets
\begin{equation}\label{Uloc}
U=\prod_{i=1}^{l}\psi_{a_{i}}(Q_{i})\times
\prod_{j=1}^{l^{\prime}}\phi_{a_{l+j}}(Q_{l+j}),
\end{equation}
both properties
(\ref{loc10}) and (\ref{loc11}) hold again. Then we have the following estimate.

\begin{lemma}\label{IPa}
Let $k,l,l^{\prime}\in \N$, $Q_{i}\in \D^{k+1,\infty
},i=1,...,l+l^{\prime}$ and set $U$ as in \eqref{Uloc}.
Consider also some $F\in(\D^{k+1,\infty})^{n}.$
Then for every $p\geq 1$ there exist some universal constants $C>0$ and $p'>p$ (depending on $k,n,p$ only) such that for every multi index $\alpha$ with
$\left\vert \alpha\right\vert \leq k$ one has
$$
\left\Vert H_{\alpha,U}(F,1)\right\Vert _{U,p}\leq C\big(1+\sigma_{U,F}%
\big(p'\big)^{k+1}\big)\Big(1+\sum_{i=1}^{l+l^{\prime}}\frac{1}{a_{i}^{k}}\left\Vert
Q_{i}\right\Vert _{k,p'}\Big)\big(1+\left\Vert F\right\Vert _{k+1,p'}^{2nk}\big).
$$

\end{lemma}

\textbf{Proof.}
For $G\in (\D^{r,p})^n$, let $|D^{(r)}G|=\sum_{\ell=0}^r\sum_{|\gamma|=\ell}|D^{\gamma}G|^2$ as in the proof of Lemma \ref{Norms1bis}. Then the following deterministic estimate for the Malliavin weights holds:
\begin{equation}\label{det-est-mall-weights}
\begin{array}{ll}
|H_{\alpha,U}(F,V)|
\leq
&\displaystyle
C\Big(\sum_{r=0}^{k}|D^{(r)}V|\Big)
\times \Big(1+\sum_{r=1}^{k}|D^{(r)}\ln U|\Big)\times\smallskip\\
&\displaystyle\times \big(1+|\det\sigma_F|^{-(k+1)}\big)\times\smallskip\\
&\displaystyle\times \Big(1+\sum_{r=1}^{k+1}|D^{(r)}F|+\sum_{r=0}^{k-1}|D^{(r)}LF|\Big)^{2nk}.
\end{array}
\end{equation}
The proof of \eqref{det-est-mall-weights} is straightforward, although non trivial, and can be found in the preprint version of the present paper, see \cite{bib:BC-Horm}. The statement now easily follows by applying to the r.h.s. of \eqref{det-est-mall-weights} the H\"older inequality and the Meyer inequality
$\|LF\|_{U, r,p}\leq \|LF\|_{ r,p}\leq C\|F\|_{r+2,p}$. $\square$

\medskip

We finally recall the result in Theorem 2.13 from \cite{[BC2]}, on which the proof of Theorem \ref{Regularity} is based.

Consider a random variable $F$, a probability measure $\Q$ and a family of probabilities $\Q_\delta$, $\delta>0$. We denote by $\mu$ the law of $F$ under $\Q$ and by $\mu_\delta$ the law of $F$ under $\Q_\delta$. In the following, we will take $\Q=\P_U$ and $\Q_\delta=\P_{U_\delta}$ as given in \eqref{loc1}, where $U$ and $U_\delta$ are both of the form \eqref{Uloc}. Actually, $\P_U$ and $\P_{U_\delta}$ are not  probability measures but they are both finite with total mass less or equal to 1, and this is enough.

We let $\E_{\Q}$ and $\E_{\Q_\delta}$ denote expectation under $\Q$ and $\Q_\delta$ respectively.

\smallskip

Fix $\delta>0$. For $m\in\N_*$ and $p\geq 1$, we say that $F\in\mathcal{R}_{m,p}(\Q_\delta)$ if
for every multi index $\alpha$ with $\left\vert \alpha\right\vert \leq m$ there exists a random variable $H_{\alpha,\delta}$ such that the following abstract integration by parts formula holds:%
\begin{equation}\label{IBP-delta}
\E_{\Q_\delta}(\partial_{\alpha}f(F))=\E_{\Q_\delta}(f(F)H_{\alpha,\delta})\qquad\forall f\in C_{c}^{\infty},\quad\mbox{with } \E_{\Q_\delta}(|H_{\alpha,\delta}|^p)<\infty.
\end{equation}
By using  Theorem 2.13 \textbf{A} in \cite{[BC2]} with $m=1$ and $k=0$, we have

\begin{theorem}\label{regularity-new}
Let $q\in \N$ and $p>1$ be fixed  and let $r_n=2(n+1)$. Let $F\in\cap_{\delta>0}\mathcal{R}_{q+3,r_n}(\Q_\delta)$.
Suppose that there exist $\theta > 0$, $C\geq 1$ and $\eta >\frac{q+n/p_{\ast }}{2}$, with $p_*$ the conjugate of $p$, such that  one has
\begin{align}
&\limsup_{\delta\to 0}
\Big(\E_{\Q_\delta}(|F|^{r_n})^{1/{r_n}}+\sum_{|\alpha|\leq q+3}\delta^{\theta|\alpha|}\E_{\Q_\delta}(|H_{\alpha,\delta}|^{r_n})^{1/{r_n}}
\Big)<\infty, \label{i13}\\
&d_{0}(\mu ,\mu _\delta )\leq C\delta ^{\eta \theta n^{2}(q+3)},
\label{i14}
\end{align}
where $d_0$ denotes  the total variation distance, that is $d_0(\mu,\nu)=\sup\{|\int fd\mu-\int fd\nu|\,:\,\|f\|_{\infty}\leq 1\}$.
Then $\mu$ is absolutely continuous and has a density $p_{F}\in W^{q,p}%
(\R^{n}).$
\end{theorem}

\textbf{Proof.}
Let us first notice that Theorem 2.13 in \cite{[BC2]} concerns a family of r.v.'s $F_\delta$, $\delta>0$, and it is assumed that all these random variables $F_\delta$ are defined on the same probability space $(\Omega, \mathcal{F}, {%
\mathbb{P}})$. But this is just for simplicity of notations. In fact the
statement concerns just the law of $(F_\delta,
H_{\alpha}(F_\delta,1),|\alpha|\leq 2m+q+1)$, where $H_{\alpha}(F_\delta,1)$ are the weights in the integration by parts formulas for $F_\delta$. So we may assume that each $%
F_\delta $ is defined on a different probability space $(\Omega_\delta,%
\mathcal{F}_\delta,{\mathbb{Q}}_\delta)$. In our case, we take $F_\delta=F$ for each $\delta$,  we work on the space $(\Omega,\mathcal{F},\Q_\delta)$ and we have $H_{\alpha}(F_\delta,1)=H_{\alpha,\delta}$.
We then apply Theorem 2.13 in \cite{[BC2]} with $m=1$ and $k=0$. (\ref{i13}) immediately gives that $\sup_\delta \E_{\Q_\delta}(|F|^{n+3})<\infty$ because $2(n+1)\geq n+3$. Moreover, in view of (2.39) in \cite{[BC2]}, the quantity $T_{q+3,2(n+1)}(F_\delta)$ in the statement of Theorem 2.13 therein can be upper bounded by
$$
S_{q+3,2(n+1)}(\delta):=
\E_{\Q_\delta}(|F|^{r_n})^{1/{r_n}}+\sum_{|\alpha|\leq q+3}\E_{\Q_\delta}(|H_{\alpha,\delta}|^{r_n})^{1/{r_n}}.
$$
As an immediate consequence of (\ref{i13}) and (\ref{i14}), all the requirements in Theorem 2.13 in \cite{[BC2]} hold, and the statement follows.
$\square$

\subsection{Proof of Theorem \ref{Regularity}}\label{proof}

We are now ready to prove our main result.

\medskip

\textbf{Proof of Theorem \ref{Regularity}.}
\textbf{Step 1: construction of the localization r.v.'s $U$ and $U_\delta$.} We consider the
functions $\psi=\psi_{1/2}$ and $\phi=\phi_2$ defined in (\ref{loc8}) and
(\ref{loc8'}) with $a=\frac{1}{2}$ and $a=2$ respectively.
We recall that in hypothesis (\ref{R3})
$ii)$ some $\gamma<\frac{1}{2}$ is considered. We denote $\lambda=\frac{1}%
{3}(\frac{1}{2}-\gamma).$ Recall that $q_{i}(W),i=1,2$ are defined in
(\ref{R5'''}). Then we define
\begin{align*}
Q_{0}  &  =r^{-1}\left\vert F-y\right\vert ,\qquad Q_{1}=\frac{68d^{3}%
}{\lambda_{\ast}\delta^{2}}G_{\delta},\qquad Q_{2}=\delta^{-(\gamma+2\lambda
)}q_{1}(W),\\
Q_{3} &=q_{2}(W), \qquad Q_{4}    =\delta^{\gamma+\lambda}\overline{a},\qquad Q_{5}=\frac
{\lambda(T,\delta)}{\lambda_{\ast}}
\end{align*}
and we set
\[
U=\psi(Q_{0}),\qquad U_{\delta}=\prod_{i=0}^{4}\psi(Q_{i})\times\phi(Q_{5}).
\]

\textbf{Step 2: construction and estimate of the weights $H_{\alpha,\delta}$ (defined in \eqref{IBP-delta}) under
$\P_{U_{\delta}}.$} We fix $k\in\N_*$ and we assume that $F\in(\D^{k+3,\infty,\infty})^n$.

Notice that for $\delta^{\lambda}\leq\frac{1}{8d}\,\sqrt{\lambda_*},$ on the set $\{U_{\delta}\neq0\}$\ we have
\[
\overline{a}(T,\delta)q_{1}(W)=\big(\delta^{\gamma+\lambda}\overline{a}(T,\delta)\big)(\delta
^{-(\gamma+2\lambda)}q_{1}(W))\times\delta^{\lambda}\leq\delta^{\lambda}%
\leq\frac{1}{8}\,\frac{\sqrt{\lambda_{\ast}}}d.
\]
The other restriction required in $\Lambda_{T,\delta}$ (see \eqref{Lambda} for the definition) are easy to check. So, we obtain
\[
\{U_{\delta}\neq0\}\subset\{\left\vert F-y\right\vert \leq r\}\cap
\Lambda_{T,\delta}.
\]
Then, by using Lemma \ref{Determinant} we have
\begin{equation}
\E_{T,\delta}\big(1_{\{U_{\delta}\neq0\}}(\det\sigma_{F,T,\delta})^{-p}\big)\leq
\frac{C_{n,p}}{\lambda_{\ast}^{np}\delta^{2np}}\label{Crit6}%
\end{equation}
where ${\sigma}_{F,T,\delta}$ is given in \eqref{R5''}.

We use the Malliavin calculus with respect to $W_{s}-W_{T-\delta}%
,s\in(T-\delta,T)$. So, we denote with $L_{\delta}$ the Ornstein  Uhlenbeck operator with respect to $W_{s}-W_{T-\delta},s\in(T-\delta,T)$ and with $\left\langle g,f\right\rangle
_{\delta}$ the scalar product in $L^{2}[T-\delta,T]$. So, ${\sigma}_{F,T,\delta}$ is the Malliavin covariance matrix of $F$ w.r.t. this partial calculus. We set, as usual, $\widehat{\sigma}_{F,T,\delta}$ the inverse of ${\sigma}_{F,T,\delta}$ and we set
$$
H_{i,U_{\delta}}(F,V):=\sum_{j=1}^{n}(V\widehat{\sigma}_{F,T,\delta}%
^{j,i}L_{\delta}F^{j}-\left\langle D(V\widehat{\sigma}_{F,T,\delta}%
^{j,i}),DF^{j}\right\rangle _{\delta}-V\widehat{\sigma}_{F,T,\delta}%
^{j,i}\left\langle D(\ln U_\delta),DF^{j}\right\rangle _{\delta}).
$$
Then \eqref{loc5} reads
$$
\E_{U_\delta}(\partial_i f(F)V)
=\E_{U_\delta}(H_{i,U_\delta}(f,V)).
$$
By iteration, for a multi index $\alpha\in\{1,\ldots,n\}^k$ we have
$$
\E_{U_\delta}(\partial_\alpha f(F)V)
=\E_{_U\delta}(H_{\alpha,U_\delta}(f,V)),
$$
where $H_{\alpha,U_\delta}(f,V)=H_{\alpha_h,U_\delta}(f,H_{(\alpha_1,\ldots,\alpha_{k-1}),U_\delta}(f,V))$. And by using Lemma \ref{IPa}, we can find $C>0$ and $p'>1$ such that
$$
\left\Vert H_{\alpha,U_\delta}(F,1)\right\Vert _{U_\delta,p}\leq C\big(1+\sigma_{U_\delta,F}%
\big(p'\big)^{k+1}\big)\Big(1+\sum_{i=1}^{5}\left\Vert
Q_{i}\right\Vert _{k,p'}\Big)\big(1+\left\Vert F\right\Vert _{k+1,p'}^{2nk}\big)
$$
with
$$
\sigma_{U_\delta,F}(p)^p
=\E_{U_\delta}((\det \sigma_{F,T,\delta})^{-p})
=\E(U_\delta(\det \sigma_{F,T,\delta})^{-p}).
$$
Since $0\leq U_\delta\leq 1_{U_\delta\neq 0}$, and by using estimate \eqref{Crit6} we get
$$
\sigma_{U_\delta,F}(p)^p
\leq \E(1_{\Lambda_{T,\delta}}(\det \sigma_{F,T,\delta})^{-p})
= \E(\E_{T,\delta}(1_{\Lambda_{T,\delta}}(\det \sigma_{F,T,\delta})^{-p}))
\leq \frac C{\lambda_*^{np}\delta^{2np}}.
$$
Moreover, by applying Remark \ref{Norms1} we obtain $\sum_{i=0}%
^{5}\left\Vert Q_{i}\right\Vert _{k+1,p^{\prime}}\leq C\delta^{-2}$. So, we
conclude that if $\left\vert \alpha\right\vert \leq k$ then
\begin{align}
\left\Vert H_{\alpha,U_{\delta}}(F,1)\right\Vert _{U_{\delta},p}
&\leq\frac
{C}{\delta^{2n(k+1)+2}}\big(1+\left\Vert F\right\Vert _{k+1,p'}^{2nk}\big)\nonumber\\
&\leq \frac{C}{\delta^{\theta k}}\big(1+\left\Vert F\right\Vert _{k+1,p'}^{2nk}\big)
\quad\mbox{with}\quad \theta=4n+2\label{I4}%
\end{align}
where $C$ is a universal constant depending on $n$, $k$ (recall that $k\geq 1$) and
$\lambda_{\ast}$.

\medskip

\textbf{Step 3: estimate of the total variation distance.} We recall that for
two non negative finite measures $\mu,\nu$ the total variation distance is
defined by
$$
d_{0}(\mu,\nu)=\sup\Big\{\Big| \int fd\mu-\int fd\nu\Big|
:\|f\|_{\infty}\leq1\Big\}.
$$
We consider the measures $\mu$ and
$\mu_{\delta}$ defined by%
\[
\int fd\mu=\E_{U}(f(F)),\qquad\int fd\mu_{\delta}=\E_{U_{\delta}}(f(F)),
\]
so that $d_0(\mu,\mu_\delta)\leq \E(|U-U_\delta|)$. Therefore, we have
\begin{align*}
d_{0}(\mu,\mu_{\delta}) \leq& \P\Big(G_{\delta}\geq\frac{\lambda_{\ast}\delta
^{2}}{68d^{3}}\Big)+\P\big(\left\vert W_{T}-W_{T-\delta}\right\vert \geq\delta
^{\frac{1}{2}-\lambda}\big)+\\
&+\P\Big(\sum_{j=1}^{d}\int_{T-\delta}^{T}\left\vert
W_{s}^{j}-W_{T-\delta}^{j}\right\vert ^{2}ds\geq\delta\Big)+\\
&
+\P\big(\overline{a}(T,\delta)   \geq\delta^{-(\gamma+\lambda)}\big)+\P\big(\{\left\vert F-y\right\vert
\leq r\}\cap\{\lambda(T,\delta)<\lambda_{\ast}\}\big)\\
=:&\sum_{i=1}^{5}%
\epsilon_{i}(\delta).
\end{align*}
For every $r\geq 1$, by using Chebychev's inequality we obtain $\epsilon
_{2}(\delta)\leq C\delta^{r(\frac{1}{2}-\gamma)}$ and in a similar
way, for every $r\geq 1$ then $\epsilon_{3}(\delta)\leq C\delta^{r/2}.$ By (\ref{R3}) $ii)$
\[
\epsilon_{4}(\delta)\leq C\delta^{r(\gamma+\lambda)}\E(\overline{a}^{r}(T,\delta))\leq
C\delta^{r\lambda}%
\]
and by (\ref{R3}) $iii)$ $\epsilon_{5}(\delta)\leq C\delta^r$ for every $r\geq 1.$ We conclude that for every $\varepsilon\geq 1$,
$$
\limsup_{\delta\to 0}\delta^{-\varepsilon}\epsilon_i(\delta)=0\quad \mbox{ for every $\varepsilon>0$ and $i=2,3,4,5$}.
$$
The behavior of $\epsilon_1(\delta)$ is given by Lemma \ref{Norms}:
if $F\in\cup_{p>6}(\D^{4,\infty,p})^n$ then there exists $\varepsilon>0$ such that $\limsup_{\delta\to 0}\delta^{-\varepsilon}\epsilon_1(\delta)=0$ and if $F\in(\D^4,\infty,\infty)$ then $\limsup_{\delta\to 0}\delta^{-\varepsilon}\epsilon_1(\delta)=0$ for every $\varepsilon>0$. Therefore,
we get
\begin{equation}\label{I5'}
\begin{array}{ll}
(i)& \mbox{$F\in\cup_{p>6}(\D^{4,\infty,p})^n$ $\Rightarrow$ $\exists$ $\varepsilon>0$ such that
$\limsup_{\delta\to 0}\delta^{-\varepsilon}d_{0}(\mu,\mu_{\delta})=0$;}\medskip\\
(ii)&\mbox{$F\in(\D^{4,\infty,\infty})^n$ $\Rightarrow$ $\forall$ $\varepsilon>0$ then
$\limsup_{\delta\to 0}\delta^{-\varepsilon}d_{0}(\mu,\mu_{\delta})=0$.}
\end{array}
\end{equation}

\medskip

\textbf{Step 4: conclusions.}
We first prove part \textbf{A} of Theorem \ref{Regularity}. Since $F\in
\cup_{p>6}(\D^{5,\infty,p})^n$, we have that  \eqref{I5'} $(i)$ holds. We apply now Theorem \ref{regularity-new} with $q=0$, $\Q=\P_U$ and $\Q_\delta=\P_{U_\delta}$. By using \eqref{I4}, \eqref{i13} holds with $\theta=4n+2$. Now, we choose $p>1$ sufficiently close to $1$ such that
$$
\Big(1-\frac 1p\Big)\times 3n^3(4n+2)<\varepsilon.
$$
So, taking $\eta=\frac{n}{p_*}$ we get $\eta>\frac{n/p_*}2$ and $3\eta\theta n^2<\varepsilon$ and by using \eqref{I5'} $(i)$ we have that
hypothesis \eqref{i14} holds. Then, by applying Theorem \ref{regularity-new}, we conclude that $\mu(dx)=f(x)dx$ and $f\in L^p(\R^n)$.

\medskip

We prove now \textbf{B} of Theorem \ref{Regularity}. As before,
\eqref{i13} holds with $\theta=4n+2$. Moreover, by
\eqref{I5'} $(ii)$, we get that \eqref{i14} holds for every choice of $p>1$ and of $\eta >\frac{q+n/p_{\ast }}{2}$. So, the only restriction in the application of Theorem \ref{regularity-new} is that $F\in\cap_{\delta>0}\mathcal{R}_{q+3,2(n+1)}(\Q_\delta)$. But in order to have this, we need that each component of $F$ is $k$-times differentiable in Malliavin sense with $k\geq(q+3)+2=q+5$, that is $q\leq k-5$. And we apply Theorem \ref{regularity-new} with $q=k-5$, giving the result.
$\square$

\section{An example from diffusion processes}\label{sect-use}

We consider the $N$ dimensional diffusion process%
\begin{equation}\label{L-diff}
dX_{t}=\sum_{j=1}^{d}\sigma_{j}(X_{t})dW_{t}^{j}+b(X_{t})dt.
\end{equation}
We assume that $\sigma_{j},b\in C_{b}^{\infty}(\R^{N}).$ In particular
$X_{T}^{i}\in\cap_{m=1}^{\infty}\D^{m,\infty,\infty}$ (see Nualart \cite{[N]}).

Our aim is to study the regularity of $\overline{X}_{T}=(X_{T}^{1}%
,...,X_{T}^{n})$ with $n\leq N.$ One may consider $\overline{X}_{t}$ as the
solution of an equation with coefficients depending on the past. We introduce
some notation. For a function $f:\R^{N}\rightarrow \R^{N}$ we denote
$\overline{f}=(f^{1},...,f^{n})$ and for $x=(x_{1},...,x_{N})\in \R^{N}$ we
denote $\overline{x}=(x_{1},...,x_{n})\in \R^{n}$ and $\widehat{x}%
=(x_{n+1},...,x_{N})\in \R^{N-n}.$ And for $\overline{x}=(x_{1},...,x_{n})\in
\R^{n}$ and $\widehat{x}=(x_{n+1},...,x_{N})\in \R^{N-n}$ we denote
$(\overline{x},\widehat{x})=(x_{1},...,x_{n},x_{n+1},...,x_{N})\in \R^{N}.$ We
define%
\begin{align*}
\Lambda_{\widehat{x},\xi}(\overline{x})  &  =\sum_{j=1}^{d}\left\langle
\overline{\sigma}_{j}(\overline{x},\widehat{x}),\xi\right\rangle ^{2}%
+\sum_{j,p=1}^{d}\left\langle \overline{[\sigma_{j},\sigma_{p}]}(\overline
{x},\widehat{x}),\xi\right\rangle ^{2}\qquad and\\
\Lambda(\overline{x})  &  =\inf_{\widehat{x}\in \R^{N-n}}\inf_{\left\vert
\xi\right\vert =1}\Lambda_{\widehat{x},\xi}(\overline{x}).
\end{align*}

\begin{proposition}
We assume that $\sigma_{j},b\in C_{b}^{\infty}(\R^{N})$ and consider a point
$\overline{x}_{0}\in \R^{n}$ such that $\Lambda(\overline{x}_0)>0.$ Then there
exists some $r>0$ such that the restriction of the law of $\overline{X}_{T}$
to $B_{r}(\overline{x}_{0})$ is absolutely continuous and has an infinitely
differentiable density on this ball.
\end{proposition}

\begin{remark}
Other types of dependence on the past may be considered. For example equations with delay (see e.g. Mohammed \cite{[Mo]}) or interacting particle systems (see e.g. L\"ocherbach \cite{eva}). For simplicity, we treat here the model given by the first $n$ components of the $N$-dimensional diffusion in \eqref{L-diff}.

\end{remark}

\textbf{Proof}. We consider $a_{j},a_{j,p},j,p=1,...,d$ defined by%
\[
a_{j}(T,\delta)=\overline{\sigma}(X_{T-\delta}),\qquad a_{j,p}(T,\delta
)=\sum_{k=1}^{N}\sigma_{j}^{k}(X_{T-\delta})\partial_{k}\overline{\sigma}%
_{p}(X_{T-\delta}).
\]
Notice that $[a]_{j,p}(T,\delta)=\overline{[\sigma_{j},\sigma_{p}%
]}(X_{T-\delta})$ so that, with the notation in (\ref{R1'}), we have
$\lambda(T,\delta)\geq$ $\Lambda(\overline{X}_{T-\delta}).$

Since the derivatives of $\sigma_{j}$ are uniformly bounded one has
$\left\vert \Lambda_{\widehat{x},\xi}(\overline{x})-\Lambda_{\widehat{x},\xi
}(\overline{x}^{\prime})\right\vert \leq C\left\vert \overline{x}-\overline
{x}^{\prime}\right\vert $ for some $C$ depending on $\left\Vert
\sigma\right\Vert _{\infty}+\left\Vert \nabla\sigma\right\Vert _{\infty}.$ So
we may find $r>0$ such that $\Lambda(\overline{x})\geq\frac{1}{2}%
\Lambda(\overline{x}_{0})$ for $\overline{x}\in B_{2r}(\overline{x}_{0}).$ It
follows that $\lambda(T,\delta)\geq\frac{1}{2}\Lambda(\overline{x}_{0})$ for
$\overline{X}_{T-\delta}\in B_{2r}(\overline{x}_{0}).$ Then
\[
\P(\{\left\vert \overline{X}_{T}-\overline{x}_{0}\right\vert <r\}\cap
\{\lambda(T,\delta)<\frac{1}{4}\Lambda(\overline{x}_{0})\})\leq \P(\left\vert
\overline{X}_{T}-\overline{X}_{T-\delta}\right\vert >r)\leq Ce^{-r^{2}%
/C^{\prime}\delta}%
\]
which proves that the hypothesis (\ref{R3}), $iii)$ holds true. Since
$\sigma_{j}$ are bounded the hypothesis (\ref{R3}), $ii)$ holds true also. Let
us check (\ref{R3}), $i)$. We compute
\[
D_{s}^{j}\overline{X}_{T}=\overline{\sigma}_{j}(X_{s})+\sum_{p=1}^{d}\int
_{s}^{T}\nabla\overline{\sigma}_{p}(X_{r})D_{s}^{j}X_{r}dW_{r}^{p}+\int
_{s}^{T}\nabla\overline{b}(X_{r})D_{s}^{j}X_{r}dr.
\]
So for $T-\delta\leq s\leq T$ we have%
\[
\E_{T,\delta}(D_{s}^{j}\overline{X}_{T})=\E_{T,\delta}(\overline{\sigma}%
_{j}(X_{s}))+\int_{s}^{T}\E_{T-\delta}(\nabla\overline{b}(X_{r})D_{s}^{j}%
X_{r})dr=a_{j}(T,\delta)+R_{\delta}^{j}(s)
\]
with%
\[
R_{\delta}^{j}(s)=\E_{T,\delta}(\overline{\sigma}_{j}(X_{s})-\overline{\sigma
}_{j}(X_{T-\delta}))+\int_{s}^{T}\E_{T,\delta}(\nabla\overline{b}(X_{r}%
)D_{s}^{j}X_{r})dr.
\]
With $L$ denoting the infinitesimal generator associated to the diffusion (\ref{L-diff}), one has
$$
\overline{\sigma}_{j}(X_{s})-\overline{\sigma
}_{j}(X_{T-\delta})
=\sum_{k=1}^d\int_{T-\delta}^s\nabla\overline{\sigma}_j(X_u)\sigma_k(X_u)dW^k_u
+\int_{T-\delta}^sL\overline{\sigma}_j(X_u)du,
$$
so that
$$
R_{\delta}^{j}(s)=
\int_{T-\delta}^s\E_{T,\delta}(L\overline{\sigma}_j(X_u))du
+\int_{s}^{T}\E_{T,\delta}(\nabla\overline{b}(X_{r}%
)D_{s}^{j}X_{r})dr.
$$
Standard computations show that $\E(|R_{\delta}^{j}(s)|^{2p})\leq C\delta^{2p}$ for any $s\in[T-\delta,T]$, so that
\begin{align*}
\E\Big(\Big|  \frac{1}{\delta}\int_{T-\delta}^{T}\left\vert \frac{\E_{T,\delta
}(D_{s}^{j}\overline{X}_{T})-a_{j}(T,\delta)}{\delta^{\frac{1}{2}+\alpha}%
}\right\vert ^{2}ds\Big|  ^{p}\Big)
&\leq\frac{1}{\delta}\int_{T-\delta}%
^{T}\E\big(\big| \delta^{-(\frac{1}{2}+\alpha)}R_{\delta}^{j}(s)\big|^{2p}\big)ds\\
&\leq C\delta^{2p(\frac{1}{2}-\alpha)}.
\end{align*}
We fix $T-\delta\leq s_{2}\leq s_{1}\leq T$ and we compute the second order
derivatives:%
\begin{align*}
\E_{T,\delta}(D_{s_{2}}^{p}D_{s_{1}}^{j}\overline{X}_{T})=\E_{T,\delta}%
(\nabla\overline{\sigma}_{j}(X_{s_{2}})D_{s_{2}}^{p}\overline{X}_{s_{1}}) &
+\sum_{k,l=1}^{d}\int_{s_{1}}^{T}\E_{T,\delta}(\partial_{k}\partial
_{l}\overline{b}(X_{r})D_{s_{2}}^{p}\overline{X}_{r}^{l}D_{s_{1}}^{j}%
\overline{X}_{r}^{k})dr\\
+\sum_{k=1}^{d}\int_{s_{1}}^{T}\E_{T,\delta}(\partial_{k}\overline{b}%
(X_{r})D_{s_{2}}^{p}D_{s_{1}}^{j}\overline{X}_{r}^{k})dr &  =a_{p,j}%
(T,\delta)+R_{\delta}^{p,j}(s_{1},s_{2})
\end{align*}
with%
\begin{align*}
R_{\delta}^{p,j} &  =\E_{T,\delta}(\nabla\overline{\sigma}_{j}(X_{s_{2}%
})D_{s_{2}}^{p}\overline{X}_{s_{1}}-\nabla\overline{\sigma}_{j}(X_{T-\delta
})\sigma(X_{T-\delta}))\\
&  +\sum_{k,l=1}^{d}\int_{s_{1}}^{T}\E_{T,\delta}(\partial_{k}\partial
_{l}\overline{b}(X_{r})D_{s_{2}}^{p}\overline{X}_{r}^{l}D_{s_{1}}^{j}%
\overline{X}_{r}^{k})dr+\sum_{k=1}^{d}\int_{s_{1}}^{T}\E_{T,\delta}%
(\partial_{k}\overline{b}(X_{r})D_{s_{2}}^{p}D_{s_{1}}^{j}\overline{X}_{r}%
^{k})dr.
\end{align*}
Similarly as before, one has $\E(|R_{\delta}^{p,j}(s_1,s_2)|^{2p})\leq C\delta^{2p}$
so that%
\begin{align*}
\frac{1}{\delta^{2}}\int_{T-\delta}^{T}&\int_{T-\delta}^{s_{1}}\E\Big(\Big|
\frac{\E_{T,\delta}(D_{s_{2}}^{p}D_{s_{1}}^{j}\overline{X}_{T})-a_{p,j}%
(T,\delta)}{\delta^{\alpha/2}}\Big|^{2p}\Big)ds_{2}ds_{1}=\\
&=\frac{1}{\delta^{2}}\int_{T-\delta}^{T}\int_{T-\delta}^{s_{1}}
\E\big(\big|\delta^{-\alpha/2}R_{\delta}^{p,j}(s_{1},s_{2})\big|^{2p}\big)ds_{2}ds_{1}\leq C\delta^{2p(1-\alpha/2)}.
\end{align*}
We conclude that for $\alpha\leq \frac{1}{2}$ we have $\varepsilon
_{\alpha,p,\delta}(a,\overline{X}_{T})\leq C$ so that the hypothesis
(\ref{R3}) $i)$ is verified. The statement now follows by applying Theorem
\ref{Regularity}. $\square$

\appendix
\section{The variance lemma}\label{app-1}
In \cite{[DY]} (see (1.f), p. 183) one gives the explicit expression of the
Laplace transform of the variance of the Brownian path on $(0,1).$ More
precisely let $B$ be an one dimensional Brownian motion and let
\begin{equation}
V(B)=\int_{0}^{1}\Big(B_{s}-\int_{0}^{1}B_{r}dr\Big)^{2}ds. \label{A1}%
\end{equation}
Then
\begin{equation}\label{A1bis}
\E(e^{-\lambda V(B)})=\frac{2\lambda}{\sinh2\lambda},\qquad\lambda>0.
\end{equation}
As an easy consequence we obtain the following estimate:

\begin{lemma}
\label{Varience}On a probability space we consider a one dimensional Brownian
motion $b$ and a random variable $r.$ We also consider two real numbers
$\alpha,\beta$ and $\delta>0$ and we denote $A_{\delta}=\{r^{2}\leq\frac
{1}{32}(\alpha^{2}+\beta^{2})\}\cap\{\left\vert \frac{1}{\delta}\int_{0}^{\delta
}b_{s}ds\right\vert \leq1\}.$ Then%
\begin{equation}
\E(1_{A_{\delta}}\exp(-\int_{0}^{\delta}(r+\alpha+\beta b_{s})^{2}ds))\leq
2\exp(-\frac{\delta^{2}}{17}(\alpha^{2}+\beta^{2})). \label{A2}%
\end{equation}

\end{lemma}

\textbf{Proof}. We consider the probability measure $\mu_{\delta}%
(ds)=\delta^{-1}1_{(0,\delta)}(s)ds$, so that
$$
\int_{0}^{\delta}(r+\alpha+\beta b_{s})^{2}ds
=\delta\int(r+\alpha+\beta b_{s})^{2}d\mu_\delta(s).
$$
Setting
\[
V_{\mu_{\delta}}(b)=\int(b_{s}-\int b_{r}d\mu_{\delta}(r))^{2}d\mu_{\delta
}(r),
\]
it is easy to check that
\begin{equation}
\int(r+\alpha+\beta b_{s})^{2}d\mu_{\delta}(s)=\left(  \int(r+\alpha+\beta
b_{s})d\mu_{\delta}(s)\right)  ^{2}+\beta^{2}V_{\mu_{\delta}}(b) \label{A3}%
\end{equation}
and
\begin{equation}
V_{\mu_{\delta}}(b)=\delta V(B)\qquad with\qquad B_{t}=\delta^{-1/2}%
b_{t\delta}. \label{A4}%
\end{equation}
We consider two cases. Suppose first that $\left\vert \alpha\right\vert
\geq4\left\vert \beta\right\vert .$ On the set $A_{\delta}$ we have
$2\left\vert \alpha\right\vert \geq\left\vert \alpha\right\vert +\left\vert
\beta\right\vert \geq8\left\vert r\right\vert $ and $\left\vert \int b_{s}%
d\mu_{\delta}(s)\right\vert \leq1$ so we obtain%
\begin{align*}
\left\vert r+\alpha+\beta\int b_{s}d\mu_{\delta}(s)\right\vert
&\geq\left\vert \alpha\right\vert -\left\vert r\right\vert -\left\vert
\beta\right\vert \left\vert \int b_{s}d\mu_{\delta}(s)\right\vert
\geq\left\vert \alpha\right\vert -\left\vert r\right\vert -\left\vert
\beta\right\vert \\
&\geq\frac{1}{2}\left\vert \alpha\right\vert \geq\frac{1}%
{4}(\left\vert \alpha\right\vert +\left\vert \beta\right\vert ).
\end{align*}

Using (\ref{A3}) this gives%
\begin{align*}
\int_{0}^{\delta}(r+\alpha+\beta b_{s})^{2}ds
&\geq \delta\left(  \int(r+\alpha+\beta b_{s})d\mu_{\delta}(s)\right)
^{2}\\
&  \geq\frac{\delta}{16}(\left\vert \alpha\right\vert +\left\vert
\beta\right\vert )^{2}
\geq\frac{\delta}{16}(\alpha^{2}+\beta^{2})\\
&\geq\frac{\delta^2}{17}(\alpha^{2}+\beta^{2}).
\end{align*}
Suppose now that $\left\vert \alpha\right\vert <4\left\vert \beta\right\vert.$
Then using (\ref{A3}) we can write
\begin{align*}
\int_{0}^{\delta}(r+\alpha+\beta b_{s})^{2}ds
&\geq \delta\beta^{2}V_{\mu_{\delta}}(b)
=\delta^2\beta^{2}V(B)\geq \frac{\delta^2}{17}(\alpha^2+\beta^2)V(B)
\end{align*}
Then we have
\begin{align*}
\E(1_{A_{\delta}}e^{-\int_{0}^{\delta}(r+\alpha+\beta b_{s})^{2}ds})
\leq &
1_{\{|\alpha|\geq 4|\beta|\}}e^{-\frac{\delta^2}{17}(\alpha^2+\beta^2)}
+1_{\{|\alpha|> 4|\beta|\}}\E(e^{-\frac{\delta^2}{17}(\alpha^2+\beta^2)V(B)})
\end{align*}
and by using (\ref{A1bis}) and the estimate $\frac{2\lambda}{\sinh(2\lambda)}\leq
2\lambda e^{-2\lambda}\leq 2e^{-\lambda}$, we get
\begin{align*}
\E(1_{A_{\delta}}e^{-\int_{0}^{\delta}(r+\alpha+\beta b_{s})^{2}ds})
\leq &
1_{\{|\alpha|\geq 4|\beta|\}}e^{-\frac{\delta^2}{17}(\alpha^2+\beta^2)}
+1_{\{|\alpha|> 4|\beta|\}}2e^{-\frac{\delta^2}{17}(\alpha^2+\beta^2)}
\end{align*}
and the statement follows.

$\square$

\section{Proof of inequality (\ref{det-est-mall-weights})}\label{app-2}
Let us briefly recall the notations we are going to use. For $r\in\N$ and a multi index $\beta\in\{1,\ldots,d\}^r$, if $F\in(\D^{r,\infty})^n$ we set
\begin{align*}
&|D^\beta F|^2=\int_{[0,T]^r}
|D^\beta_{s_1\ldots s_r}F|^2ds_1\cdots ds_r
=\sum_{j=1}^d\int_{[0,T]^r}|D^\beta_{s_1\ldots s_r}F^j|^2ds_1\cdots ds_r
\mbox{ and }\\
&|D^{(r)} F|^2=\sum_{|\beta|=r}|D^\beta F|^2.
\end{align*}
For the sake of completeness, we allow $\beta=\emptyset$, or equivalently $|\beta|=0$: we  set
$$
D^\beta F=F\quad\mbox{and}\quad|D^{(0)} F|^2=|F|^2.
$$
Moreover, $\<\cdot,\cdot\>$ denotes the scalar product in $L^2([0,T],dt)$, so that for $F,G\in(\D^{1,\infty})^n$,
$$
\<DG,DF\>=\sum_{i=1}^d\int_{[0,T]}D^i_sGD^i_sFds
=\sum_{i=1}^d\sum_{j=1}^n\int_{[0,T]}D^i_sG^jD^i_sF^jds
$$
and $|DF|^2=\<DF,DF\>=|D^{(1)}F|^2$.

For $F$ taking values in $\R^n$, $V$ in $\R$ and $\alpha$ multi index of length $k$ in $\{1,\ldots,n\}$, let $H_{\alpha,U}(F,V)$ denote the weight in \eqref{loc7}, that is the weight from the integration by parts formula of order $k$ of $F$ w.r.t. $V$ localized through $U$. The appendix is devoted to the proof of the following
\begin{proposition}\label{p1}
For $\ell=0,1,\ldots$, let $\beta\in\{1,\ldots,d\}^\ell$ (the case $\ell=0$ referring to $\beta=\emptyset$) and for $k=1,2,\ldots$, let $\alpha\in\{1,\ldots,n\}^k$ be a multi index of length $k$. On the set $\{U>0\}$, let $V,\ln U\in\D^{k+\ell,\infty}$ and let $F\in(\D^{k+\ell+1,\infty})^n$ be such that the associated Malliavin covariance matrix $\sigma_F$ is invertible. Then, on the set $\{U>0\}$ the following estimate holds:
\begin{align*}
|D^\beta H_{\alpha,U}(F,V)|
\leq
&C\Big(\sum_{r=0}^{k+\ell}|D^{(r)}V|\Big)
\times \Big(1+\sum_{r=1}^{k+\ell}|D^{(r)}\ln U|\Big)\times\\
&\times \big(1+|\det\sigma_F|^{-(k+\ell+1)}\big)\times\\
&\times \Big(1+\sum_{r=1}^{k+\ell+1}|D^{(r)}F|+\sum_{r=0}^{k+\ell-1}|D^{(r)}LF|\Big)^{2(k+\ell)n},
\end{align*}
$C$ being a positive constant depending on $\beta$ and $\alpha$ but independent of $U$, $F$  and $V$.
\end{proposition}

As a consequence, taking $\beta=\emptyset$ one gets that  \eqref{det-est-mall-weights} holds.
The proof of Proposition \ref{p1} requires some preliminary estimates.

\begin{lemma}\label{l1}
Let $r\in\N$ and $\gamma\in\{1,\ldots,d\}^r$ be a multi index of length $r$. Then for every $F,G\in (\D^{r+1,\infty})^n$ the following statements hold:
\begin{align}
&|D^\gamma\<DG,DF\>|\leq 2\Big(\sum_{\ell=1}^{r+1}|D^{(\ell)}G|\Big)
\Big(\sum_{\ell=1}^{r+1}|D^{(\ell)}F|\Big),\label{l1.1}\\
&|D^\gamma\sigma_F|\leq c
\Big(\sum_{\ell=1}^{r+1}|D^{(\ell)}F|\Big)^{2}\label{l1.2}\\
&|D^\gamma\det\sigma_F|\leq c_r
\Big(\sum_{\ell=1}^{r+1}|D^{(\ell)}F|\Big)^{2n}\label{l1.3}\\
&|D^\gamma(\det\sigma_F)^{-1}|\leq c_r
\Big(1+|\det\sigma_F|^{-(r+1)}\Big)\Big(1+\sum_{\ell=1}^{r+1}|D^{(\ell)}F|\Big)^{2nr}\label{l1.4}\\
&|D^\gamma\widehat\sigma_F|\leq c_r
\Big(1+|\det\sigma_F|^{-(r+1)}\Big)\Big(1+\sum_{\ell=1}^{r+1}|D^{(\ell)}F|\Big)^{2n(r+1)-2}\label{l1.5}
\end{align}
Here, $c$ and $c_r$ denote suitable positive constants, possibly depending on $r$ but universal w.r.t. the choice of $F$ and/or $G$.
\end{lemma}

\textbf{Proof.} \emph{Proof of \eqref{l1.1}.} One has
\begin{align*}
D^\gamma_{s_1,\ldots, s_r}\<DG,DF\>
&=\<D^\gamma_{s_1,\ldots, s_r}D_\cdot G, D_\cdot F\>+\<D_\cdot G, D^\gamma_{s_1,\ldots, s_r}D_\cdot F\>
\end{align*}
so that by the Cauchy-Schwartz inequality we get
\begin{align*}
|D^\gamma_{s_1,\ldots, s_r}\<DG,DF\>|
&\leq |D^\gamma_{s_1,\ldots, s_r}D_\cdot G||DF| +|D G| |D^\gamma_{s_1,\ldots, s_r}D_\cdot F|.
\end{align*}
By noticing that $|D^\gamma_{s_1,\ldots, s_r}D_\cdot G|^2=
\int_{[0,T]}|D^\gamma_{s_1,\ldots, s_r}D_s G|^2ds$, the statement  follows.

\smallskip

\emph{Proof of \eqref{l1.2}.} Since $|D^\gamma\sigma_F^{ij}|
=|D^\gamma\<DF^i,DF^j\>|$, the result follows from \eqref{l1.1}.

\smallskip

\emph{Proof of \eqref{l1.3}.} Recall that
$\det \sigma_F=\sum_{\rho\in\mathcal{P}_n}\sigma_F^{1\rho_1}\cdots\sigma_F^{n\rho_n},$
where $\mathcal{P}_n$ is the set of all permutations of $(1,\ldots,n)$. For $\gamma$ multi index in $\{1,\ldots,d\}$ with $|\gamma|=r$ we set $s_\gamma\in\R^r$ as $s_\gamma=(s_{\gamma_1},\ldots,s_{\gamma_r})$. Then, we can write
$$
D^\gamma_{s_\gamma}\det \sigma_F
=\sum_{\rho\in\mathcal{P}_n}D^\gamma\big(\sigma_F^{1\rho_1}\cdots\sigma_F^{n\rho_n}\big)
=\sum_{\rho\in\mathcal{P}_n}
\sum_{\beta_1,\ldots,\beta_n\in\mathcal{A}_\gamma}
D^{\beta_1}_{s_{\beta_1}}\sigma_F^{1\rho_1}\cdots D^{\beta_n}_{s_{\beta_n}}\sigma_F^{n\rho_n}
$$
where ``$\beta_1,\ldots,\beta_n\in\mathcal{A}_\gamma$'' means that $\beta_1,\ldots,\beta_n$ is a partition of $\gamma$
running through the list of all of the ``blocks'' of $\gamma$.  We use now \eqref{l1.3} and we obtain
$$
|D^\gamma\det \sigma_F|
\leq
\sum_{\rho\in\mathcal{P}_n}
\sum_{\beta_1,\ldots,\beta_n\in\mathcal{A}_\gamma}
|D^{\beta_1}\sigma_F^{1\rho_1}|\cdots |D^{\beta_n}\sigma_F^{n\rho_n}|
\leq c \Big(\sum_{\ell=1}^{r+1}|D^{(\ell)}F|\Big)^{2n}.
$$
\emph{Proof of \eqref{l1.4}.} We set again $s_\gamma=(s_{\gamma_1},\ldots,s_{\gamma_r})$. For $f\in C^r$ we can write
$$
D^\gamma_{s_\gamma}f(G)
=\sum_{\ell=1}^rf^{(\ell)}(G)\sum_{\beta_1,\ldots,\beta_\ell\in\mathcal{B}_\gamma}
D^{\beta_1}_{s_{\beta_1}}G\cdots D^{\beta_\ell}_{s_{\beta_\ell}}G
$$
where ``$\beta_1,\ldots,\beta_\ell\in\mathcal{B}_\gamma$'' means that $\beta_1,\ldots,\beta_\ell$ are non empty multi indexes of $\gamma$ running through the list of all of the (non empty) ``blocks'' of $\gamma$. Then, it follows that
$$
|D^\gamma f(G)|
\leq
\sum_{\ell=1}^r|f^{(\ell)}(G)|\!\!\!\sum_{
\beta_1,\ldots,\beta_\ell\in\mathcal{B}_\gamma}
\!\!|D^{\beta_1}G|\cdots |D^{\beta_\ell}G|
\leq c\max_{1\leq \ell\leq r}
|f^{(\ell)}(G)|\,\Big(1+\sum_{j=1}^r |D^{(j)}G|\Big)^{r}
$$
because $1\leq |\beta_i|\leq |\gamma|=r$. We consider now $f(x)=1/x$ and $G=\det\sigma_F$. Here, $|f^{(\ell)}(x)|= \ell! x^{-(\ell+1)}$. So, by noticing that for $\ell\leq r$
$$
1+|\det\sigma_F|^{-(\ell+1)}\leq 2(1\vee |\det\sigma_F|^{-1})^{\ell+1}\leq 2(1\vee |\det\sigma_F|^{-1})^{r+1}
\leq 2(1 + |\det\sigma_F|^{-(r+1)})
$$
and by using \eqref{l1.3} one gets the result.

\smallskip

\emph{Proof of \eqref{l1.5}.}
We set $\tilde \sigma_F$ as the matrix of cofactors, so that $\widehat\sigma_F^{ij}=(-1)^{i+j}(\det\sigma)^{-1}\tilde\sigma_F^{ji}$. Then,
$$
D^\gamma_{s_{\gamma}}\widehat\sigma_F^{ij}=(-1)^{i+j}\sum_{\beta_1,\beta_2\in \mathcal{A}_\gamma}
D^{\beta_1}_{s_{\beta_1}}(\det\sigma_F^{-1})D^{\beta_2}_{s_{\beta_2}}\tilde\sigma_F^{ji}
$$
where we say that $\beta_1,\beta_2\in\mathcal{A_\gamma}$ iff $\beta_1,\beta_2$ is a partition of $\gamma$.  By recalling that $\tilde\sigma_F^{ji}$ is the determinant of the sub-matrix of $\sigma_F$ obtaining by deleting the $j$th row and the $i$th column of $\sigma$, we can apply \eqref{l1.3} to $D^{\beta_2}\tilde\sigma_F^{ji}$. And by using \eqref{l1.4} for $D^{\beta_1}(\det\sigma_F^{-1})$, \eqref{l1.5} immediately holds. $\square$

\bigskip

We are now ready for the

\medskip

\textbf{Proof of Proposition \ref{p1}.}
We first consider the case $|\alpha|=1$. Here, we use a reduced notation and we write
\begin{align*}
H_{U}(F,V)
&=V\widehat \sigma_{F}LF-\widehat\sigma_F\<DV,DF\>-V\<D\widehat\sigma_F,DF\>
-V\widehat\sigma_{F}\< D\ln U,DF\>
\end{align*}
(recall that $V$ is always one dimensional, while $F$ takes values in $\R^n$), so that
$H_{i,U}(F,V)$ is the $i$th entry of the random vector $H_{U}(F,V)$.
For a multi index $\beta$ with $|\beta|=\ell$ we have
\begin{align*}
D^\beta_{s_\beta} H_{U}(F,V)
=&
\sum_{\gamma_1,\gamma_2,\gamma_3\in \mathcal{A}^3_\beta}D^{\gamma_1}_{s_{\gamma_1}}VD^{\gamma_2}_{s_{\gamma_2}}\widehat \sigma_{F}D^{\gamma_3}_{s_{\gamma_3}}LF
-\sum_{\gamma_1,\gamma_2\in \mathcal{A}^2_\beta}D^{\gamma_1}_{s_{\gamma_1}}\widehat\sigma_FD^{\gamma_2}_{s_{\gamma_2}}\<DV,DF\>+\\
&-\sum_{\gamma_1,\gamma_2\in \mathcal{A}^2_\beta}D^{\gamma_1}_{s_{\gamma_1}}VD^{\gamma_2}_{s_{\gamma_2}}\<D\widehat\sigma_F,DF\>+\\
&-\sum_{\gamma_1,\gamma_2\gamma_3\in \mathcal{A}^3_\beta}D^{\gamma_1}_{s_{\gamma_1}}VD^{\gamma_2}_{s_{\gamma_2}}\widehat\sigma_{F}D^{\gamma_3}_{s_{\gamma_3}}\< D\ln U,DF\>
\end{align*}
where the condition ``$\gamma_1,\ldots,\gamma_i\in\mathcal{A}^i_\beta$'' means that
$\gamma_1,\ldots,\gamma_i$ is a partition of $\beta$ given by $i$ subsets.
We set now
$$
\mathcal{H}_k(G)=\sum_{r=0}^k|D^{(r)}G|.
$$
Then, for a suitable constant $C$ (independent of $V$, $F$ and $U$) that can vary from line to line, we can write
\begin{align*}
|D^\beta H_{U}(F,V)|
\leq&
C\Big(\mathcal{H}_\ell(V)\mathcal{H}_\ell(\widehat\sigma_F)\mathcal{H}_\ell(LF)
+\mathcal{H}_\ell(\widehat\sigma_F)\mathcal{H}_\ell(\<DV,DF\>)+\\
&+\mathcal{H}_\ell(V)\mathcal{H}_\ell(\<D\widehat\sigma_F,DF\>)+\mathcal{H}_\ell(V)\mathcal{H}_\ell(\widehat\sigma_{F})\mathcal{H}_k(\< D\ln U,DF\>)\Big)
\end{align*}

We estimate the above terms by using Lemma \ref{l1}:

\smallskip

$\bullet$ from \eqref{l1.5} one has
$$
\mathcal{H}_\ell(\widehat\sigma_F)
\leq C
\big(1+|\det\sigma_F|^{-(\ell+1)}\big)\Big(1+\sum_{i=1}^{\ell+1}|D^{(i)}F|\Big)^{2n(\ell+1)};
$$
$\bullet$
from \eqref{l1.1} one has
$$
\mathcal{H}_\ell(\<DV,DF\>)
\leq C\sum_{i=1}^{\ell+1}|D^{(i)}V|\times
\sum_{i=1}^{\ell+1}|D^{(i)}F|
\leq
C\sum_{i=1}^{\ell+1}|D^{(i)}V|\times
\Big(1+\sum_{i=1}^{\ell+1}|D^{(i)}F|\Big);
$$
$\bullet$
from \eqref{l1.1} and \eqref{l1.5} one has
\begin{align*}
\mathcal{H}_\ell(\<D\widehat\sigma_F,DF\>)
&\leq C\sum_{i=1}^{\ell+1}|D^{(i)}\widehat\sigma_F|\times
\sum_{i=1}^{\ell+1}|D^{(i)}F|\\
&\leq C\big(1+|\det\sigma_F|^{-(\ell+2)}\big)
\Big(1+\sum_{i=1}^{\ell+2}|D^{(i)}F|\Big)^{2n(\ell+1)};
\end{align*}
$\bullet$
from \eqref{l1.1} one has
$$
\mathcal{H}_\ell(\<D\ln U,DF\>)
\leq C\sum_{i=1}^{\ell+1}|D^{(i)}\ln U|\times
\sum_{i=1}^{\ell+1}|D^{(i)}F|.
$$
So, by inserting the above estimates we get the result for $|\alpha|=1$. The case $|\alpha|=k>1$ now easily follows by induction. $\square$

\end{document}